\documentclass{article}
\usepackage[utf8]{inputenc}
\usepackage[margin=1in]{geometry}

\usepackage[table, svgnames, dvipsnames]{xcolor}

\usepackage[utf8]{inputenc} 
\usepackage[T1]{fontenc}    
\usepackage{hyperref}                
\usepackage{booktabs}       
\usepackage{amsfonts}       
\usepackage{nicefrac}       
\usepackage{microtype}      
\usepackage{float}
\usepackage{pagecolor}
\usepackage{amsmath, amssymb, amsthm, mathtools}
\usepackage[authoryear]{natbib}
\usepackage{amsmath}

\emergencystretch=\maxdimen
\numberwithin{equation}{section}

\hypersetup{
    colorlinks=true,
    linkcolor=blue,
    filecolor=magenta,      
    urlcolor=blue,
}

\newcommand{\EX}{\mathbb{E}}
\DeclareMathOperator*{\argmin}{arg\,min}

\newcommand{\doublenabla}{%
  \nabla\mkern-11mu\nabla%
}

\renewcommand{\text}[1]{\textnormal{#1}}

\allowdisplaybreaks

\title{Complex normalizing flows can almost be information Kähler-Ricci flows}

\vspace{0.4in}

\author{%
  Andrew Gracyk\thanks{\small Department of Mathematics, Purdue University, West Lafayette, IN 47907, United States, \texttt{agracyk@purdue.edu}}
}
\date{}

\begin{document}

\maketitle

\begin{abstract}
\noindent We develop interconnections between the complex normalizing flow for data drawn from Borel probability measures on the twofold realification of the complex manifold and a nonlinear flow nearly Kähler-Ricci. The complex normalizing flow relates the initial and target realified densities under the complex change of variables, necessitating the log determinant of the ensemble of Wirtinger Jacobians. The Ricci curvature of a Kähler manifold is the second order mixed Wirtinger partial derivative of the log of the local density of the volume form. Therefore, we reconcile these two facts by drawing forth the connection that the log determinant used in the complex normalizing flow matches a Ricci curvature term under differentiation and conditions. The log density under the normalizing flow is kindred to a spatial information metric under an augmented Jacobian and a Bayesian perspective to the parameter, thus under the continuum limit the log likelihood matches a Fisher metric, or more closely a Kähler cross-entropy Hessian. This recovers a Kähler-Ricci flow variation up to a time derivative and expectation, or an average-valued Kähler-Einstein flow. Using this framework, we establish other relevant results, attempting to bridge the statistical and ordinary behaviors of the complex normalizing flow to the geometric features of our derived Kähler flow.
\end{abstract}

\medskip

\vspace{2mm}

\noindent \textbf{Key words.} Information geometry, complex geometry, Kähler-Ricci flow, Kähler geometry, normalizing flow, Fisher information, information metric, log determinant, Ricci curvature, complex diffeomorphism, holomorphic, anti-holomorphic, Bayesian, instantaneous change of variables, Mabuchi, Perelman, surgery

\vspace{2mm}

\noindent\textbf{AMS MSC Classifications (2020):}  	53B35, 53B12, 53Z50, 34M04, 34M45

\tableofcontents

\section{Introduction}

We contribute to the geometric modality for the normalizing flow with an unusual connection to Kähler geometry via the complex Ricci flow. Complex Ricci curvature has a distinction from geometric structure with open sets diffeomorphic to Euclidean space through a closed form using log determinants of the metric \cite{song2012lecturenoteskahlerricciflow}. The normalizing flow is recognizable in machine learning literature for its use of the density transformation law, utilizing the log determinant of the Jacobian in its log likelihood evolution. Both complex Ricci curvature and a change of variables transformation law share in concord use of such a log determinant, thus our work resides in drawing this connection and reconciling the consistency among these two frameworks via this log determinant.

\vspace{2mm}

\noindent The overarching aim of the normalizing flow is to learn a pushforward of data drawn from (realified) Borel probability measures \cite{inbook} with finite second moments, the initial data set easily sampleable and the target, possibly nontrivial to estimate or derive in a parametric regression task, having known and given data. These objectives are compatible with generative modeling \cite{lipman2023flowmatchinggenerativemodeling} \cite{papamakarios2021normalizingflowsprobabilisticmodeling} \cite{zhai2025normalizingflowscapablegenerative} \cite{grathwohl2018ffjordfreeformcontinuousdynamics}, since by learning a closed form map from a simple to mature density, we can sample the simple density to generate a new sample in the mature density space. Normalizing flows have connections to literature outside of generative modeling, for example they have connections to neural ordinary differential equations \cite{chen2019neuralordinarydifferentialequations} \cite{scagliotti2024normalizingflowsapproximationsoptimal} \cite{xu2022infinitelydeepbayesianneural}, mean field games \cite{Huang_2023} \cite{zhang2023meanfieldgameslaboratorygenerative}, variational inference \cite{rezende2016variationalinferencenormalizingflows}, anomaly detection \cite{rosenhahn2024quantumnormalizingflowsanomaly}, Bayesian statistics \cite{Roch_2026}, and representation learning. Normalizing flows are primarily established for measures on real-valued data, but we offer new perspective to that which is complex-valued.

\vspace{2mm}

\noindent Our primary results will be done via the information-theoretic perspective. Is is well established that the pushforward flow map of a traditional normalizing flow is a diffeomorphism at each increment \cite{brehmer2020flowssimultaneousmanifoldlearning} \cite{ross2021tractabledensityestimationlearned}, thus the data at each time can implicitly be treated as those of a manifold. In section \ref{sec:normalizing_flow_geometry}, we establish the manners in which we consider the geometry of the normalizing flow. We can geometrically detail this diffeomorphism, or biholomorphism in the complex case, via the metric pullback. Therefore, the log determinant of the transformation law uses Riemannian, or Hermitian, metric information. This metric coincides with a Fisher metric through use of the log likelihood, but it is standard in Fisher information theory that the Fisher metric is taken with respect to the parameter and the data is annihilated as an argument via the  integration. We invert how we treat the manifold, and we integrate over the parameter space instead and differentiate with respect to data. Thus our Fisher arguments somewhat deviate from those often found to be more standard, and we transpose the manner of the information manifolds.

\vspace{2mm}

\noindent In order to offer an information-theoretic perspective via an inverted Fisher information, we provide a Bayesian perspective to the parameter of the pushforward neural networks. To integrate-out the parameter in the Fisher metric, we treat the parameter as a distributional quantity \cite{xu2022infinitelydeepbayesianneural}
\cite{yamauchi2023normalizingflowsbayesianposteriors} \cite{trippe2018conditionaldensityestimationbayesian} \cite{maroñas2021transforminggaussianprocessesnormalizing}. Thus, we assume the parameter follows a posterior measure which is conditioned on the total dataset. Moreover, we will assume the posterior is affected by the time of the normalizing flow \cite{NEURIPS2022_e9e1a0ab} \cite{NEURIPS2020_df1a336b}, for example consistent with applications in online Bayesian learning or continual learning, which is more of an assumption for generalized purposes. Our approaches are easily consistent when the posterior is set fixed in time, thus our approaches are generalized Bayesian statistics. We advance the Bayesian paradigm for the normalizing flow consequently.

\vspace{2mm}

\noindent We examine to various levels the curvature uniformization property and how that manifests in the normalizing flow. To introduce, we crucially remark that our flow is more reminiscent of a Kähler-Einstein condition, and that we do not claim the normalizing flow is a Kähler-Ricci flow exactly. It is well known from Ricci flow theory that nontrivial manifolds diffuse to those of constant curvature \cite{chen2005noteuniformizationriemannsurfaces} \cite{chau2007surveykahlerricciflowyaus} \cite{chen2002uniformizationtheoremcompletenoncompact} \cite{chau2003gradientkahlerriccisolitonsuniformization}, although this result is more elaborate than we present it here because this is also affected by manifold dimension, surgical qualities, etc. \cite{BamlerRicciFlowNotes2015} \cite{Topping2006}. In the normalizing flow, this has some manifestation via the exchange of a complicated density to that which is sampleable and potentially isotropic, i.e. a unit Gaussian, or a flat metric. It can be noted curvature is created in the generative task $q_0 \rightarrow q_{K}$, thus it is notable to keep a sign convention in mind. We remark this curvature property in the normalizing is rather by convention than necessity since the initial density is by choice: there is no intrinsic requirement $q_0$ is simple or flat.

\vspace{2mm}

\noindent Throughout this work, we will establish several connections to two conventions of the normalizing flow, being the baseline/discrete complex normalizing flow \cite{tran2019discreteflowsinvertiblegenerative} and the complex continuous normalizing flow \cite{grathwohl2018ffjordfreeformcontinuousdynamics} \cite{kingma2018glowgenerativeflowinvertible} \cite{dinh2017densityestimationusingreal}. The continuous version adopts the instantaneous change of variables theorem \cite{chen2019neuralordinarydifferentialequations}, and does several other things for us, which are: (1) simplifying computations via a continuity equation; (2) and allows use of such a theorem in our proofs, which means that, under a time-varying manifold,
\begin{align}
\frac{d}{dt} \log p(z(t)) = - \text{div}_{\omega_t}(f) - \text{Tr}_{\omega_t} \left( \frac{\partial \omega_t}{\partial t} \right)
\end{align}
for suitable $f$. We will employ these two techniques in our proofs, thus we highlight we consider both versions of these normalizing flows.

\vspace{2mm}

\noindent Lastly, we note that standard real-valued normalizing flows can be extended to complex normalizing flows in certain cases via the connection that $x+\sqrt{-1}y$ can be treated as real numbers $(x,y)$, i.e. $\mathbb{C}^d$ is realified and isomorphic to $\mathbb{R}^{2d}$. Thus our techniques can partially be applied to typical normalizing flows as well, as long as $d$ is even. Thus, we remark the manifold as in $\ref{fig:complexnf_and_curvature}$ is a 2-real manifold but a 1-complex manifold when orientable.

\section{Notations and conventions}
\label{sec:notations}

We will denote $\overline{z}$ to be the complex conjugate of $z$ and $\dagger$ to be a Hermitian transpose. We will denote $\Psi$ to be a map in the normalizing flow, and $\Phi \in  \{ \phi \in C^\infty(\mathbb{C}^d \times [0,T]; \mathbb{R}) : \sqrt{-1}\partial\overline{\partial}\phi(\cdot,t) > 0, \int_{\mathbb{C}^d} |\phi|^2 e^{-\phi} \frac{(\sqrt{-1}\partial\overline{\partial}\phi)^d}{d!} < \infty \ \forall t \}$ to be a Kähler potential. We use
\begin{align}
\partial \overline{\partial} \Phi = \sum_{ij} \frac{\partial \overline{\partial} \Phi}{\partial z^i \partial \overline{z}^j} dz^i \wedge d\overline{z}^j 
\end{align}
as the form version with Dolbeault operators $(\partial, \overline{\partial})$. We will denote
\begin{align}
\frac{\partial}{\partial z^j} = \frac{1}{2} \Big( \frac{\partial}{\partial x^j} - \sqrt{-1} \frac{\partial}{\partial y^j} \Big), \frac{\partial}{\partial \overline{z}^j} = \frac{1}{2} \Big( \frac{\partial}{\partial x^j} + \sqrt{-1}\frac{\partial}{\partial y^j} \Big)
\end{align}
to be the Wirtinger derivatives. A critical assumption that we will make is a Bayesian one on our neural network parameter $\theta$, and that $\theta \sim p(\theta|\mathcal{D},t)$ for posterior $p$. We will assume every Hermitian metric $h(t) \in \Gamma(M, T^{*(1,0)} M \otimes T^{*(0,1)}M)$ is Kähler based on the formulation of \ref{sec:normalizing_flow_geometry}. We will not pull back biholomorphic functions \cite{Greb_2019} \cite{marini2022crrelativeskaehlermanifolds}, since the log determinant of a holomorphic pullback has vanishing Ricci curvature in equivalent dimensions. In particular, pushforwards of isotropic metrics with holomorphic functions do not admit nondegenerate Ricci curvature (see Appendix \ref{app:vanishing_ricci} for proof). Here, $T^{*(1,0)} M$ denotes the holomorphic cotangent bundle. We will denote $\omega$ the Kähler form representative of the de Rham cohomology class $[\omega]$, and so $d\omega=0$, $d: \Omega^k(M) \to \Omega^{k+1}(M)$. We will use without providing the proof the instantaneous change of variables theorem as in \cite{chen2019neuralordinarydifferentialequations} so that it holds in the manifold complex case, but we will also use the real-argument Euclidean version. We will use the manifold divergence \cite{kreutzdelgado2009complexgradientoperatorcrcalculus} in local coordinates
\begin{align}
\text{div}_{\omega_t}(\varphi) = \frac{1}{\det(h)} \left( \sum_j \frac{\partial}{\partial z^j} \big( \det(h) \varphi^j \big) + \sum_j \frac{\partial}{\partial \overline{z}^j} \big( \det(h) \overline{\varphi^j} \big) \right) ,
\end{align}
under the assumption $\varphi$ is a smooth (1,0)-vector field. When $\varphi$ is holomorphic in the vector-valued sense, the above simplifies. Observe this definition is consistent with the Kähler Laplacian, but the $\det(h)$ is absorbed in the Laplacian definition. We will use the $\sqrt{-1}$ over the $i$ imaginary unit convention to distinguish from index notation. We will use notation $\Phi = \Phi_t$ to denote time-dependence. We will use the integral conventions 
\begin{align}
\int_M \frac{\omega^d}{d!} = \int_M \det h \left( \frac{\sqrt{-1}}{2} \right)^d dz_1 \wedge d\overline{z}_1 \wedge \dots \wedge dz_d \wedge d\overline{z}_d = \int_M d \text{Vol}_h 
\end{align}
since $\omega^d = d! \left( \frac{\sqrt{-1}}{2} \right)^d \bigwedge_i \zeta^i \wedge \overline{\zeta}^i$ with respect to an orthonormal frame, and $\int_{\mathbb{C}^d} f(z) d\mu(z) = \int_{\mathbb{R}^{2d}} f(x,y) d\mu(x,y)$, where $z_k = x_k + \sqrt{-1} y_k$.

\section{Background}

\subsection{Normalizing flows}

We will restrict our analysis to discrete normalizing flows for now, since discrete flows use a Jacobian transformation law, which will be our means of Ricci curvature. A discrete complex normalizing flow seeks a target data distribution through a series of non-holomorphic complex diffeomorphisms \cite{rezende2016variationalinferencenormalizingflows}
\begin{align}
\label{eqn:composition}
z_K = \Psi_{K,\theta} \circ \Psi_{K-1,\theta} \circ \hdots \circ \Psi_{1,\theta} (z_0) .
\end{align}
$z_i \in \mathbb{C}^d$,  $\Psi_{i,\theta} : M_i \subseteq \mathbb{C}^{d} \rightarrow M_{i+1} \subseteq \mathbb{C}^{d}$ is a non-holomorphic smooth diffeomorphism between open sets, although we will slightly relax the open set condition variously throughout this work. The underlying data distribution evolves according to the change of variables
\begin{align}
\label{eqn:transformation_law_normalizing_flow}
\log q_{K,\theta}(z_K) = \log q_{0,\theta}(z_0) - \sum_{k=1}^K \log \Big| \det \mathcal{J}_{k,\theta} \Big| .
\end{align}
Here, we have the augmented Jacobian
\begin{align}
\mathcal{J} = \begin{pmatrix}
\nabla_z \Psi & \nabla_{\overline{z}} \Psi \\
\nabla_z \overline{\Psi} & \nabla_{\overline{z}} \overline{\Psi}
\end{pmatrix} = \begin{pmatrix}
\frac{\partial \Psi^1}{\partial z^1} & \cdots & \frac{\partial \Psi^1}{\partial z^d} & \frac{\partial \Psi^1}{\partial \overline{z}^1} & \cdots & \frac{\partial \Psi^1}{\partial \overline{z}^d} \\
\vdots & \ddots & \vdots & \vdots & \ddots & \vdots \\
\frac{\partial \Psi^d}{\partial z^1} & \cdots & \frac{\partial \Psi^d}{\partial z^d} & \frac{\partial \Psi^d}{\partial \overline{z}^1} & \cdots & \frac{\partial \Psi^d}{\partial \overline{z}^d} \\
\frac{\partial \overline{\Psi}^1}{\partial z^1} & \cdots & \frac{\partial \overline{\Psi}^1}{\partial z^d} & \frac{\partial \overline{\Psi}^1}{\partial \overline{z}^1} & \cdots & \frac{\partial \overline{\Psi}^1}{\partial \overline{z}^d} \\
\vdots & \ddots & \vdots & \vdots & \ddots & \vdots \\
\frac{\partial \overline{\Psi}^d}{\partial z^1} & \cdots & \frac{\partial \overline{\Psi}^d}{\partial z^d} & \frac{\partial \overline{\Psi}^d}{\partial \overline{z}^1} & \cdots & \frac{\partial \overline{\Psi}^d}{\partial \overline{z}^d}
\end{pmatrix} .
\end{align}
The above simplifies under a holomorphic map, which we do not restrict here for our reasons as in \ref{sec:notations}. We prove in Appendix \ref{app:change_of_variables}
\begin{align}
|\det \mathcal{J}_k| = \frac{\det h_{k-1}}{\det h_{k}} .  
\end{align}

\subsection{The Kähler-Ricci flow}

The Kähler-Ricci flow is the evolution equation \cite{song2012lecturenoteskahlerricciflow} (up to constants)
\begin{align}
\frac{\partial}{\partial t} h_{i \overline{j}} = -\text{Ric}_{i \overline{j}}(h) \ \ \ \ \ \text{or equivalently} \ \ \ \ \ \frac{\partial}{\partial t} \omega = - \text{Ric}(\omega) ,
\end{align}
where $\omega = \omega(t)  \in \Gamma(M, \Lambda^{(1,1)} T^*M)$ is a (1,1)-form family of Kähler metrics with Hermitian metric (1,1)-tensor representation $h$, and $\text{Ric}$ is the (complex) Ricci curvature which satisfies as a (1,1)-tensor
\begin{align}
\label{eqn:ricci_curvature}
\text{Ric}_{i \overline{j}} = - \partial_i \partial_{\overline{j}} \log \det(h) .
\end{align}
Again, $\partial_i \partial_{\overline{j}}$ are the mixed Wirtinger partial derivatives. It can be noted the Ricci form corresponding to $\omega$ uses $\rho = - \sqrt{-1} \partial \overline{\partial} \log \det (h)= \sqrt{-1} \sum_{ij} \text{Ric}_{i \overline{j}} dz^i \wedge d \overline{z}^j$. We say $h$ is Kähler if it satisfies $\partial_k h_{i \overline{j}} = \partial_i h_{k \overline{j}}$, or equivalently if $h_{i \overline{j}} = \partial_i \partial_{\overline{j}} \Phi$ for suitable potential $\Phi$.

\subsection{The geometry of the normalizing flow}
\label{sec:normalizing_flow_geometry}

We reconcile the complex normalizing flow and its geometry. We treat the measure in which the data derives from the volume form that follows $q_t = e^{-\varphi_t}$. To ensure this is a valid distribution, we normalize via Boltzmann-type
\begin{align}
q_t(z) = \frac{1}{Z_t} e^{-\varphi_t(z)} ,  Z_t = \int_{M} e^{-\varphi(z, \overline{z},t)} \left( \frac{\sqrt{-1}}{2} \right)^d \bigwedge_{i} dz^i \wedge d\overline{z}^i  = \int_{M_{2\mathbb{R}}} e^{-\varphi(x, y,t)} dx^1 \wedge dy^1 \wedge \dots \wedge dx^d \wedge dy^d .
\end{align}
We decouple potential $\varphi_t$ from the Kähler potential $ \Phi \in  \{ \phi \in C^\infty(\mathbb{C}^d \times [0,T]; \mathbb{R}) : \sqrt{-1}\partial\overline{\partial}\phi(\cdot,t) > 0, \int_{\mathbb{C}^d} |\phi|^2 e^{-\phi} \frac{(\sqrt{-1}\partial\overline{\partial}\phi)^d}{d!} < \infty \ \forall t \}$ and $\Phi \in \text{SPSH}(\mathbb{C}^d)$. $\Phi$ is defined such that
\begin{align}
h_{i \overline{j}} = \partial_i \partial_{\overline{j}} \Phi .
\end{align}
We enforce the relationships $q \propto \det h$ and $\dot{\Phi} = \log ( \frac{q}{p} )$, where $p$ is a target density. The latter equation is a parabolic Monge-Ampère flow given a target. It can be noted Kähler-Ricci flow is deterministic based on initial data, and does not use a target $p$, while $\dot{\Phi} = \log ( \frac{q}{p} )$ describes evolution towards a target and is a consequence of the normalizing flow.

\vspace{2mm}

\noindent In practice, $h_{i \overline{j}} = \partial_i \partial_{\overline{j}} \Phi$ is trivially true for an information metric by taking $\Phi = \EX[ -\log q ]$, but this is not guaranteed to be globally plurisubharmonic, which allows $h$ to be Kähler, so instead we assume $h$ is close to Kähler, or at least locally. It can be noted $\Phi \in \text{SPSH}(\mathbb{C}^d)$ when $q$ is unit Gaussian. For this work, we will work with $h$ to be locally Kähler and not globally, and so we will assume it holds for specific regions of our work and not everywhere. For example, for two moons dataset, $h$ is locally Kähler on nondegenerate portions of data but not globally. We also remark not all of our results necessitate $h$ to be Kähler, but it is true all of our results essentially require $h$ to be derivative of a potential. Also, we remark our information metric is defined so that the expectation is taken with respect to a different density than in its interior. By the above criteria, it is still a valid metric, although this is a notable distinction and it is somewhat unconventional.

\section{Our primary contribution}
\label{sec:main_contribution}

\begin{figure}[htbp]
  \centering
  \includegraphics[width=1.0\textwidth]{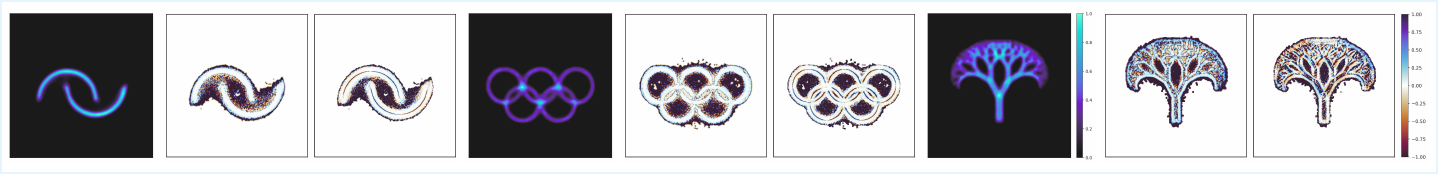}
  \vspace{-2mm}
  \caption{We plot results from a complex normalizing flow on the complex (1) two moons; (2) Olympic rings; (3) fractal tree datasets on the left with the complex density using $(\text{Re}(z_K),\text{Im}(z_K))$. On the right, we plot Kähler scalar curvature $R= -2 (h_{z \overline{z}})^{-1} \partial_{z \overline{z}} \log h_{z \overline{z}}$ and holomorphic score curvature proxy $\widetilde{R} = - \frac{1}{p} ( \partial_{xy} \log p ) / ( | \nabla \log p|^2 + \epsilon)$. We have normalized both scales to $[0.0,1.0]$, and we restrict the extreme values to the 60th percentile to prevent disproportionality. We use a $\sigma=1.0$ parameter to control smoothing of the histogram before computing scalar curvature on the right.}
  \label{fig:complexnf_and_curvature}
\end{figure}

\noindent Our primary contribution is drawing the connection between the complex Ricci curvature of \ref{eqn:ricci_curvature} and the normalizing flow discretized with \ref{eqn:transformation_law_normalizing_flow}. Since $h$ is positive definite, it is true that $\log |\det(h)| = \log \det(h)$. Identifying the probability density with the metric volume form $q_{k,\theta} = \det h_{k,\theta}$, and taking the complex Hessian, we see
\begin{align}
\label{eqn:transformation_law_Kähler}
\partial_i \partial_{\overline{j}} \left( \log q_{K,\theta}(z) - \log q_{0,\theta}(z) \right) &= -\partial_i \partial_{\overline{j}} \sum_{k=1}^K \log |\det \mathcal{J}_{k,\theta}| = \text{Ric}_{i\overline{j}}(h_{0,\theta}) - \text{Ric}_{i\overline{j}}(h_{K,\theta}),
\end{align}
which is a pointwise identity, where $h_{k,\theta}$ is a suitable Hermitian/Kähler metric at iterate $k$. $\mathcal{J}$ is the augmented Jacobian, and $h$ is the information metric. It is crucial in the above that the map $\Psi$ is not a biholomorphism, as this will lead to a splitting of the log determinant into holomorphic and anti-holomorphic parts. As a consequence, this creates vanishing Ricci curvature.

\vspace{2mm}

\noindent Let us simplify the left-hand side of \ref{eqn:transformation_law_Kähler}. Let us examine
\begin{align}
&\mathbb{E}_{\theta \sim p(\theta|\mathcal{D},t)} \left[ \partial_i \partial_{\overline{j}} \left( \log q_{K,\theta}(z) - \log q_{0,\theta}(z) \right) \right] = \mathbb{E}_{\theta \sim p(\theta|\mathcal{D},t)} [\partial_i \partial_{\overline{j}} \log q_{K,\theta}(z)] - \mathbb{E}_{\theta \sim p(\theta|\mathcal{D},t)} [\partial_i \partial_{\overline{j}} \log q_{0,\theta}(z)].
\end{align}
We can note this is equal to
\begin{align}
-I_{i\overline{j},K}(z,t) + I_{i\overline{j},0}(z,t) = \mathbb{E}_{\theta \sim p(\theta|\mathcal{D},t)} [\partial_i \partial_{\overline{j}} \log q_{K,\theta}(z)] - \mathbb{E}_{\theta \sim p(\theta|\mathcal{D},t)} [\partial_i \partial_{\overline{j}} \log q_{0,\theta}(z)],
\end{align}
which acts as a spatial information Hermitian metric. Substituting the sum from the right-hand side of \ref{eqn:transformation_law_Kähler}, we conclude
\begin{align}
-I_{i\overline{j},K}(z,t) + I_{i\overline{j},0}(z,t) = \mathbb{E}_{\theta \sim p(\theta|\mathcal{D},t)} [\text{Ric}_{i\overline{j}}(h_{0,\theta}) - \text{Ric}_{i\overline{j}}(h_{K,\theta})],
\end{align}
\noindent or equivalently over a single iterate $k$,
\begin{align}
I_{i\overline{j},k}(z,t) - I_{i\overline{j},k-1}(z,t) = \mathbb{E}_{\theta \sim p(\theta|\mathcal{D},t)} [\text{Ric}_{i\overline{j}}(h_{k,\theta}) - \text{Ric}_{i\overline{j}}(h_{k-1,\theta})].
\end{align}
In the continuum limit, dividing by $\Delta t$ and taking $\Delta t \to 0^+$, this describes the rate of change of the information metric
\begin{align}
\lim_{\Delta t \to 0^+} \frac{I_{i\overline{j},k} - I_{i\overline{j},k-1}}{\Delta t} = \lim_{\Delta t \rightarrow 0^+} \frac{\mathbb{E}_{\theta \sim p(\theta|\mathcal{D},t)} [\text{Ric}_{i\overline{j}}(h_{k,\theta}) - \text{Ric}_{i\overline{j}}(h_{k-1,\theta})]}{\Delta t}.
\end{align}
Noticing that $h$ is our information metric, the change of variables recovers the flow
\begin{align}
\partial_t h_{i\overline{j}} = \partial_t \mathbb{E}_{\theta \sim p(\theta|\mathcal{D},t)}  [ \text{Ric}_{i\overline{j}}(h_{t,\theta}) ].
\end{align}
This demonstrates that under the standard map, the metric evolves in concord with the expected time derivative of its Ricci curvature. Alternatively, this result implies in the time-independent case
\begin{align}
h_{i\overline{j}}(z, t) = \mathbb{E}_{\theta \sim p(\theta|\mathcal{D})} [\text{Ric}_{i\overline{j}}(h_{t,\theta})] + C_{i\overline{j}}(z) .
\end{align}
This is only valid when the posterior is not dependent on time. It is crucial that the expectation is here, which causes averaging. Without it, we have a type of Kähler-Einstein condition. We refer to \ref{app:kahler_einstein_conditions} for more discussion on Kähler-Einstein conditions. Moreover, this equation, were it to hold pointwise, would be nontraditional. The expectation ensemble average saves this result from being anomalous, since this pseudo-flow would force the metric to be almost Kähler-Einstein at every iterate.

\vspace{2mm}

\noindent For literature relating pullbacks to Fisher metrics, we reference \cite{holbrook2017nonparametricfishergeometrychisquare, ay2017parametrizedmeasuremodels, Itoh_2023, Bruveris_2018, Facchi_2010, cho2025statisticalbergmangeometry}. We have mostly bypassed the use of the traditional pullback via the augmented Jacobian $J$, since the Ricci curvature vanishes under a holomorphic pullback. Thus, $\mathcal{J}$ acts similarly to pullback replacement.

\vspace{2mm}

\noindent The framework we just established is not as easily extended in the continuous case because the continuous case obeys the instantaneous change of variables and implicitly uses Wirtinger Jacobian and its variations but not explicitly. Instead, our argument that the above holds in the continuous case resides in the fact that a discrete normalizing flow matches a continuous normalizing flow in the continuum limit \cite{chen2019neuralordinarydifferentialequations} \cite{salman2018deepdiffeomorphicnormalizingflows}.

\vspace{2mm}

\noindent To reiterate, the expectation is taken with respect to $p$ but the interior uses $q$. Thus, this is unconventional in Fisher metric literature; however, by what we established earlier, it is still a valid metric by the mixed Wirtinger derivative and the plurisubharmonic qualities. For example, the concept of a metric being compatible with a cross-entropy Hessian
\begin{align}
H(\theta) = \mathbb{E}_{x \sim p_{\text{data}}} [-\nabla^2 \log q_\theta(x)]
\end{align}
is not entirely new in literature \cite{shoham2025flatnessall}, thus our techniques are reconciled with these existing Fisher-Rao alternatives.

\subsection{Training}
\label{sec:training}

In this section, we outline the training procedure of the (complex) normalizing flow and draw its connections to our results.

\vspace{2mm}

\noindent Let $z \sim p_{\alpha}$ be a sample drawn from a base density in complex space, i.e. $z \in \mathbb{C}^d$. Let $\Psi_\theta$ be a non-holomorphic neural network corresponding to the totality of the composition as in \ref{eqn:composition}. The transformation law via a change of variables is given by
\begin{align}
q(w) = \frac{ p_\alpha (\Psi_\theta^{-1}(w)) }{ | \det \mathcal{J} _{\theta} | } .
\end{align}
Taking the log,
\begin{align}
\label{eqn:log_density_transformation}
\log q(w) = \log p_\alpha (\Psi_\theta^{-1}(w)) - \log | \det \mathcal{J}_{\theta} |_{z = \Psi_{\theta}^{-1}(w)} | .
\end{align}
The optimization objective is given by
\begin{align}
\argmin_\theta \text{KL} (p_\beta \parallel \Psi_\theta \# p_\alpha) = \argmin_\theta - \EX_{w \sim p_\beta} \Big[ \log ( \Psi_\theta \# p_\alpha(w)) \Big] .
\end{align}
The equality follows after ignoring constant terms. Substituting in \ref{eqn:log_density_transformation} as recall the pushforward of density formula is a change of variables, the loss is
\begin{align}
 \argmin_\theta - \EX_{w \sim p_\beta} \Big[ \log p_\alpha (\Psi_\theta^{-1}(w)) - \log | \det \mathcal{J}_{\theta} |_{z = \Psi_{\theta}^{-1}(w)} | \Big] .
\end{align}
Assuming the base distribution is standard normal complex Gaussian $p_\alpha(z) = \pi^{-d} \text{exp}(- z^{\dagger} z )$, and ignoring constants, we arrive at
\begin{align}
& \argmin_\theta \EX_{w \sim p_\beta} \Big[ | \Psi_\theta^{-1} |^2 + \log | \det \mathcal{J}_{\theta} |_{z = \Psi_{\theta}^{-1}(w)} | \Big] 
\\
& = \argmin_\theta \int_{M} \left( | \Psi_\theta^{-1} |^2 + \log | \det \mathcal{J}_{\theta} |_{z = \Psi_{\theta}^{-1}(w)} | \right)  p_{\beta}(w)  ( \frac{\sqrt{-1}}{2} )^d \bigwedge_i dw^i \wedge d\overline{w}^i .
\end{align}
Differentiating inside the objective,
\begin{align}
& \EX_{w \sim p_\beta} \Big[ \partial_i \partial_{\overline{j}} \Big( | \Psi_\theta^{-1} |^2 + \log | \det \mathcal{J}_{\theta} |_{z = \Psi_{\theta}^{-1}(w)} | \Big) \Big] 
\\
= \ & \EX_{w \sim p_\beta} \Big[ \partial_i \partial_{\overline{j}}  | \Psi_\theta^{-1} |^2 +  \partial_i \partial_{\overline{j}} \Delta \log  \det_{\mathbb{C}} (h)   \Big] =  \EX_{w \sim p_\beta} \Big[ \partial_i \partial_{\overline{j}}  | \Psi_\theta^{-1} |^2 -  \Delta \text{Ric}_{i \overline{j}}   \Big]  .
\end{align}
We have used $\Delta$ as shorthand notation for a difference formula in the discrete case. It can be noted dominated convergence is not applicable in the above.

\section{Additional theoretical results}

\begin{figure}[h]
  \makebox[\textwidth][c]{\includegraphics[width=0.8\textwidth, page=1]{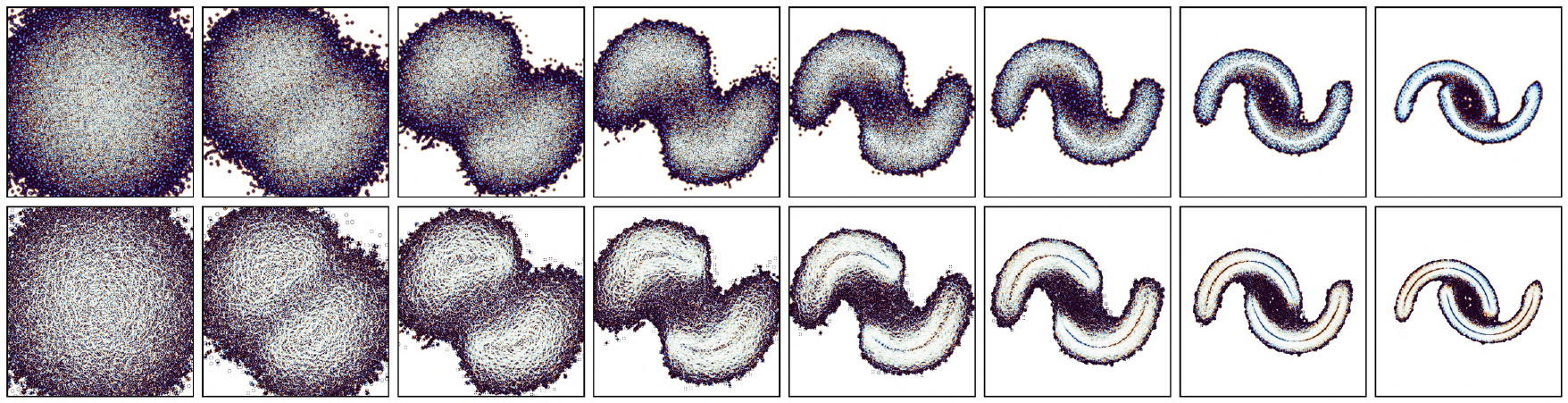}}\\[0.1em]
  \makebox[\textwidth][c]{\includegraphics[width=0.8\textwidth, page=2]{complexnf_layers_continuous_curvature.pdf}}\\[0.1em]
  \makebox[\textwidth][c]{\includegraphics[width=0.8\textwidth, page=3]{complexnf_layers_continuous_curvature.pdf}}
  \vspace{-2mm}
  \caption{We plot timesteps along the complexified continuous normalizing flow with the curvature quantities corresponding to Figure \ref{fig:complexnf_and_curvature}.}
  \label{fig:complexnf_layers_comparison_curvature}
\end{figure}

\noindent \textbf{Theorem 1.} \textit{ Denote $h(t) \in \Gamma(M, T^{*(1,0)} M \otimes T^{*(0,1)}M)$ the Kähler information metric, $\Phi \in  \{ \phi \in C^\infty(\mathbb{C}^d \times [0,T]; \mathbb{R}) : \sqrt{-1}\partial\overline{\partial}\phi(\cdot,t) > 0, \int_{\mathbb{C}^d} |\phi|^2 e^{-\phi} \frac{(\sqrt{-1}\partial\overline{\partial}\phi)^d}{d!} < \infty \ \forall t \}$ the Kähler potential. Suppose $\frac{\partial}{\partial t} \Psi_{\theta} \in \Gamma^{\infty}(T \mathbb{C}^d) $ is as in a complex continuous normalizing flow. Suppose our derived Kähler-Einstein flow is satisfied. Suppose $\theta \sim p(\theta_0|\mathcal{D},t)$ follows a posterior, and that $p,q$ follow the complex instantaneous change of variables with respect to $f,g$. Then we have, up to expectation on the posterior, the following:
\begin{itemize}
\item Assume local coordinates are time-independent. Then Ricci curvature statistically obeys
\begin{align}
\text{Ric}_{i \overline{j}} = - \EX_p \Big[ \text{div}_{\mathbb{R}^n}(f)(\partial_i \partial_{\overline{j}} \log q ) \Big] + \EX_p \Big[ \partial_i \partial_{\overline{j}} ( \text{div}_{\omega_t}(g) + \partial_t \log \omega_t) \Big] ,
\end{align}
and moreover, the Ricci curvature time derivative obeys
\begin{align}
\frac{\partial}{\partial t} \text{Ric}_{\ell \overline{k}} = -\partial_\ell \partial_{\overline{k}} \Bigg( h^{\overline{j} i} \Big( \EX_p \Big[ \text{div}_{\mathbb{R}^n}(f)(\partial_i \partial_{\overline{j}} \log q ) \Big] + \EX_p \Big[ \partial_i \partial_{\overline{j}} ( \text{div}_{\omega_t}(g) + \partial_t \log \omega_t) \Big] \Big) \Bigg) .
\end{align}
We have noted $p$ is independent of the manifold.
\item The time derivative of scalar curvature obeys
\begin{align}
\frac{\partial}{\partial t} R  & = - h^{\overline{j} \ell} h^{\overline{k} i} \partial_t \text{Ric}_{\ell \overline{k}} \text{Ric}_{i \overline{j}} + h^{\overline{j} i} \partial_i \partial_{\overline{j}} ( h^{\overline{\ell} k} \text{Ric}_{k \overline{\ell}} ) ,
\end{align}
which is consistent with \cite{song2012lecturenoteskahlerricciflow}.
\item Under the (complex) instantaneous change of variables theorem, the particle vector field obeys
\begin{align}
V = - h^{\overline{j} i} \partial_{\overline{j}} \dot{\Phi} \frac{\partial}{\partial z^i}  ,
\end{align}
up to constants.
\item Denote $p$ the terminal density and $q_t$ the density induced by $\Psi$ at a corresponding step. The first derivative of the KL divergence obeys
\begin{align}
\frac{d}{d t} \text{KL} ( q_t \parallel p) &= - \int_M | \partial \dot{\Phi} |^2 q_t \frac{\omega_t^d}{d!} + \int_M \dot{q}_t \frac{\omega_t^d}{d!} = - \text{Fisher information} + \text{correction} .
\end{align}
Denote $p$ the terminal density and $q_t$ the density induced by $\Psi$ at a corresponding step. The first derivative of the KL divergence obeys
\begin{align}
\frac{d}{d t} \text{KL} ( q_t \parallel p) &= - \int_M | \partial \dot{\Phi} |^2 q_t \frac{\omega_t^d}{d!} + \int_M \dot{q}_t \frac{\omega_t^d}{d!} = - \text{Fisher information} + \text{correction} .
\end{align}
Suppose (1) assume $p$ is log-concave $-\sqrt{-1} \partial \overline{\partial} \log p \geq \lambda \omega$ for some $\lambda > 0$; (2) assume $\sqrt{-1} \partial \overline{\partial} \log \frac{q_t}{p} \geq 0$; (3) and the manifold is closed. The second derivative obeys
\begin{align}
\frac{d^2}{dt^2}\text{KL}(q_t\|p) \geq & \left( 2 \lambda   - \frac{C}{8 \nu} \right) \int_M|\partial \dot{\Phi}|^2 q_t \frac{\omega_t^d}{d!}
\\
&  + 2 \text{Re} \int_M \langle \partial \ddot{\Phi}, \partial \dot{\Phi} \rangle_{\omega_t} q_t \frac{\omega_t^d}{d!}
\\ 
& + \int_M \text{Ric}(\partial \dot{\Phi}, \overline{\partial} \dot{\Phi} ) q_t \frac{\omega_t^d}{d!}
\\
& - ( 1 + 2 \nu) \int_M |\nabla^{2,0} \dot{\Phi}|^2 q_t \frac{\omega_t^d}{d!} 
\\
& + \int_M  (\Delta_{\overline{\partial}} \dot{\Phi}) \langle \partial \dot{\Phi}, \partial \log q_t \rangle_{\omega_t} q_t \frac{\omega_t^d}{d!} 
\\
& + \int_M \Big( -(\text{Tr}_{\omega_t}(\mathcal{R}))^2 q_t + h^{\overline{j} \ell} h^{\overline{k} i} \mathcal{R}_{\ell \overline{k}} \mathcal{R}_{i \overline{j}} q_t + h^{\overline{j} i} \partial_t \mathcal{R}_{i \overline{j}} q_t + \text{Tr}_{\omega_t}(\mathcal{R}) \dot{q}_t \Big) \frac{\omega_t^d}{d!} ,
\end{align}
where $C,\nu$ are real constants independent of dimension, and $\mathcal{R}$ is suitably defined. 
\end{itemize}
}

\noindent \textbf{Theorem 2.} \textit{Denote $h(t) \in \Gamma(M, T^{*(1,0)} M \otimes T^{*(0,1)}M)$ the Kähler information metric, $\Phi \in  \{ \phi \in C^\infty(\mathbb{C}^d \times [0,T]; \mathbb{R}) : \sqrt{-1}\partial\overline{\partial}\phi(\cdot,t) > 0, \int_{\mathbb{C}^d} |\phi|^2 e^{-\phi} \frac{(\sqrt{-1}\partial\overline{\partial}\phi)^d}{d!} < \infty \ \forall t \}$ the Kähler potential. Suppose $\Psi_{\theta}$ is as in either a complex continuous or discrete normalizing flow. Suppose our derived Kähler-Einstein flow is satisfied. Suppose $\theta \sim p(\theta_0|\mathcal{D},t)$ follows a posterior. Then we have, up to expectation on the posterior, the following:
\begin{itemize}
\item The density evolution under our condition but normalized obeys \begin{align}
\log q_{k,\theta} = \log q_{k-1,\theta} - \Big[ \log |\det(\mathcal{J}_{k,\theta})| + \nu \Delta t \log q_{k-1,\theta} \Big]  .
\end{align}
\item If $\mathcal{M}$ is the Mabuchi functional, then it has the same critical points as the KL divergence at a Kähler scalar curvature condition.
\item Under Kähler-Ricci flow, suitable $f$ and the Dirichlet metric, $\Phi$ is governed by a gradient flow satisfying
\begin{align}
\dot{\Phi} = -\doublenabla \mathcal{F}_K = \log \det h - f  .
\end{align}
Here, $\doublenabla$ is the gradient w.r.t. the Dirichlet metric. Similar results are found in \cite{Cao1985} \cite{phong2007multiplieridealsheaveskahlerricci} \cite{Huang_2020} \cite{chen2009stabilitykahlerricciflow}.
\end{itemize}
}

\section{Experiments}

\begin{figure}[htbp]
  \vspace{0mm}
  \centering
  \includegraphics[scale=0.45]{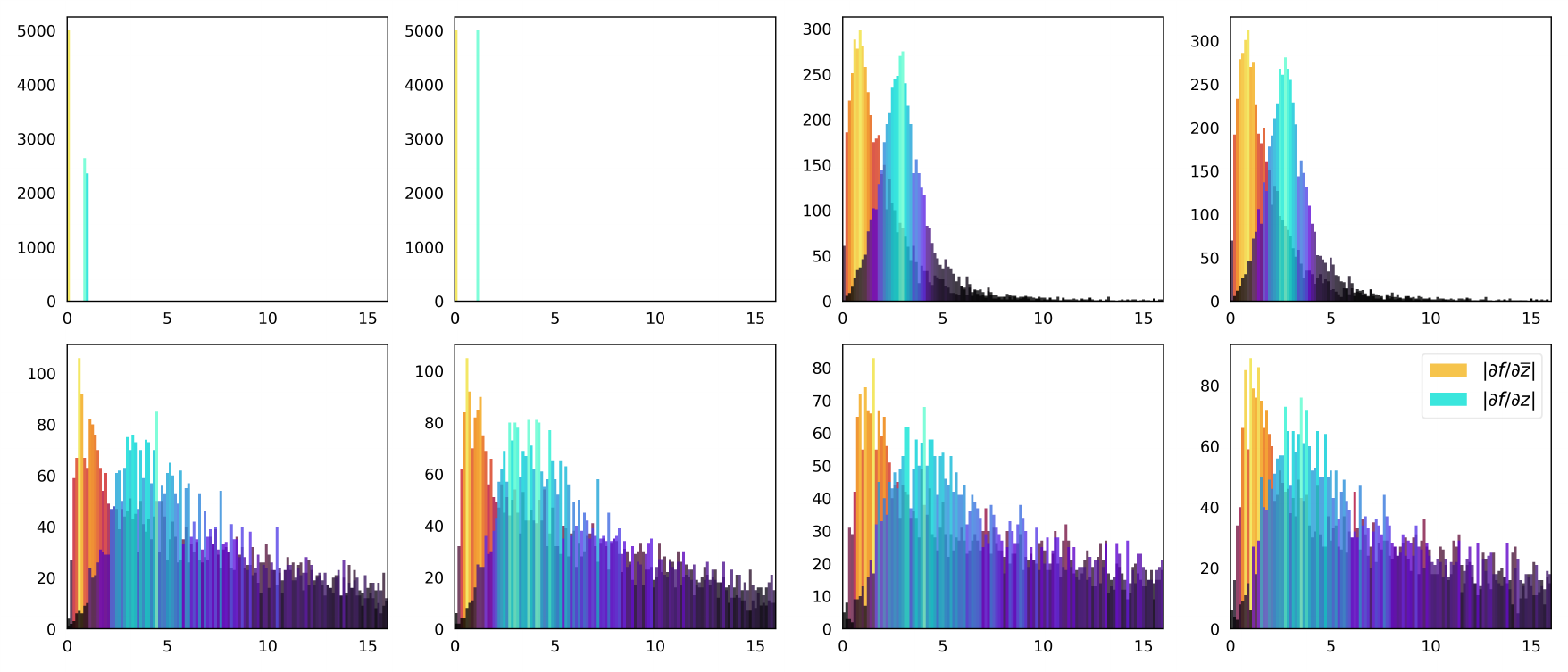}
  \caption{We illustrate a holomorphic condition or a lack thereof on each layer $\Psi_{k,\theta}$ of the 8-layer complex discrete normalizing flow after sampling 5,000 points according to complex unit Gaussian (the starting distribution in the complex normalizing flow). We plot on the fractal tree dataset. We desire $|\partial f / \partial \overline{z}|$ (oranges) to be closer to zero than $|\partial f / \partial z|$ (blues). A holomorphic condition is adverse due to vanishing Ricci curvature, thus it is desirable orange is not zero exactly.}
  \label{fig:holomorphic_layer_plot_fractaltree}
\end{figure}

\noindent \textbf{Architectures for flows.} We discuss our architecture as in Figure \ref{fig:complexnf_and_curvature}. We use a complex GELU activation $\text{GELU}(z) \cdot z /|z|$ which is not holomorphic, and we found best success with this activation as opposed to others in experimentation, for example softplus analogs. The complex linear layer implement a layer in the neural network forward as
\begin{align}
(A + \sqrt{-1}B)(z_r + \sqrt{-1} z_i) = (Az_r - Bz_i) + \sqrt{-1}(Az_i + Bz_r) ,
\end{align} 
where $A,B$ are weights. We use a coupling layer with complex affine transformation after a layer $z_1' = z_1 \cdot e^{s(z_{0})} + t(z_{0}), s,t : \mathbb{C} \rightarrow \mathbb{C}, z_1 \in \mathbb{C}$ and the exponential is the complex exponential. The Jacobian calculation only uses $2 \cdot \text{Re}(s)$ since we require
\begin{align}
J = 
\begin{pmatrix} 
I & 0 \\ 
\frac{\partial}{\partial z_0} \left( z_1 e^{s(z_0)} + t(z_0) \right) & e^{s(z_0)} 
\end{pmatrix} \in \mathbb{M}^{4 \times 4} ,
\end{align}
thus the determinant only requires the diagonal and the condition $\det_{\mathbb{R}} (J) = |\det_{\mathbb{C}} J|^2$ is not crucial with this architecture. In particular, $(z_0, z_1)_{\text{out}} = \Psi((z_0, z_1)_{\text{in}})$, i.e. each layer $\Psi_{k,\theta} : \mathbb{C}^2 \rightarrow \mathbb{C}^2$.

\section{Conclusions and limitations}

We developed connections between the complex normalizing flow and the Kähler-Ricci flow. We showcased that under Wirtinger differentiation and conditions, such as a Bayesian neural network parameter, the two mathematical frameworks hold equivalences. Our methods are primarily interesting because Kähler-Ricci flow is nonstandard in machine learning literature. Our methods establish that we can construct normalizing flows through a complex geometric lens, thus we have combined two research areas that are dissimilar. Thus, results typically found in Kähler literature can be applied variously to normalizing flows upon suitable circumstance. Our primary limitation is that, while our mathematical frameworks develop connections between two disjoint research areas, there is not much space to develop computational results as consequences of our methods.

\section{Acknowledgments}

I gratefully acknowledge financial support from Purdue University Department of Mathematics. I would like to thank Rongjie Lai at Purdue University in the Department of Mathematics for his reading course collaboration and helpful discussions on normalizing flows. I would also like to thank Nicholas McCleerey at Purdue University in the Department of Mathematics for his graduate differential geometry course, which overall helped with this work.

\begingroup

\bibliographystyle{plainnat}
\bibliography{bibliography}

@misc{rezende2016variationalinferencenormalizingflows,
      title={Variational Inference with Normalizing Flows}, 
      author={Danilo Jimenez Rezende and Shakir Mohamed},
      year={2016},
      eprint={1505.05770},
      archivePrefix={arXiv},
      primaryClass={stat.ML},
      url={https://arxiv.org/abs/1505.05770}, 
}

@misc{song2012lecturenoteskahlerricciflow,
      title={Lecture notes on the K\"ahler-Ricci flow}, 
      author={Jian Song and Ben Weinkove},
      year={2012},
      eprint={1212.3653},
      archivePrefix={arXiv},
      primaryClass={math.DG},
      url={https://arxiv.org/abs/1212.3653}, 
}

@misc{zhang2018mongeampereflowgenerativemodeling,
      title={Monge-Amp\`ere Flow for Generative Modeling}, 
      author={Linfeng Zhang and Weinan E and Lei Wang},
      year={2018},
      eprint={1809.10188},
      archivePrefix={arXiv},
      primaryClass={cs.LG},
      url={https://arxiv.org/abs/1809.10188}, 
}

@misc{xu2022infinitelydeepbayesianneural,
      title={Infinitely Deep Bayesian Neural Networks with Stochastic Differential Equations}, 
      author={Winnie Xu and Ricky T. Q. Chen and Xuechen Li and David Duvenaud},
      year={2022},
      eprint={2102.06559},
      archivePrefix={arXiv},
      primaryClass={stat.ML},
      url={https://arxiv.org/abs/2102.06559}, 
}

@misc{yamauchi2023normalizingflowsbayesianposteriors,
      title={Normalizing Flows for Bayesian Posteriors: Reproducibility and Deployment}, 
      author={Yukari Yamauchi and Landon Buskirk and Pablo Giuliani and Kyle Godbey},
      year={2023},
      eprint={2310.04635},
      archivePrefix={arXiv},
      primaryClass={nucl-th},
      url={https://arxiv.org/abs/2310.04635}, 
}

@misc{trippe2018conditionaldensityestimationbayesian,
      title={Conditional Density Estimation with Bayesian Normalising Flows}, 
      author={Brian L Trippe and Richard E Turner},
      year={2018},
      eprint={1802.04908},
      archivePrefix={arXiv},
      primaryClass={stat.ML},
      url={https://arxiv.org/abs/1802.04908}, 
}

@misc{maroñas2021transforminggaussianprocessesnormalizing,
      title={Transforming Gaussian Processes With Normalizing Flows}, 
      author={Juan Maroñas and Oliver Hamelijnck and Jeremias Knoblauch and Theodoros Damoulas},
      year={2021},
      eprint={2011.01596},
      archivePrefix={arXiv},
      primaryClass={cs.LG},
      url={https://arxiv.org/abs/2011.01596}, 
}

@article{Greb_2019,
   title={Canonical complex extensions of Kähler manifolds},
   volume={101},
   ISSN={1469-7750},
   url={http://dx.doi.org/10.1112/jlms.12287},
   DOI={10.1112/jlms.12287},
   number={2},
   journal={Journal of the London Mathematical Society},
   publisher={Wiley},
   author={Greb, Daniel and Wong, Michael Lennox},
   year={2019},
   month=nov, pages={786–827} }

@misc{marini2022crrelativeskaehlermanifolds,
      title={CR relatives Kaehler manifolds}, 
      author={Stefano Marini and Michela Zedda},
      year={2022},
      eprint={2109.13039},
      archivePrefix={arXiv},
      primaryClass={math.DG},
      url={https://arxiv.org/abs/2109.13039}, 
}

@misc{george2025complexmongeampereequationpositive,
      title={Complex Monge-Amp\`ere equation for positive $(p,p)$-forms on compact K\"ahler manifolds}, 
      author={Mathew George},
      year={2025},
      eprint={2411.06497},
      archivePrefix={arXiv},
      primaryClass={math.AP},
      url={https://arxiv.org/abs/2411.06497}, 
}

@misc{chu2024kahlermanifoldsnonnegativemixed,
      title={On K\"ahler manifolds with non-negative mixed curvature}, 
      author={Jianchun Chu and Man-Chun Lee and Jintian Zhu},
      year={2024},
      eprint={2408.14043},
      archivePrefix={arXiv},
      primaryClass={math.DG},
      url={https://arxiv.org/abs/2408.14043}, 
}

@misc{chen2019neuralordinarydifferentialequations,
      title={Neural Ordinary Differential Equations}, 
      author={Ricky T. Q. Chen and Yulia Rubanova and Jesse Bettencourt and David Duvenaud},
      year={2019},
      eprint={1806.07366},
      archivePrefix={arXiv},
      primaryClass={cs.LG},
      url={https://arxiv.org/abs/1806.07366}, 
}

@misc{holbrook2017nonparametricfishergeometrychisquare,
      title={The nonparametric Fisher geometry and the chi-square process density prior}, 
      author={Andrew Holbrook and Shiwei Lan and Jeffrey Streets and Babak Shahbaba},
      year={2017},
      eprint={1707.03117},
      archivePrefix={arXiv},
      primaryClass={stat.ME},
      url={https://arxiv.org/abs/1707.03117}, 
}

@misc{ay2017parametrizedmeasuremodels,
      title={Parametrized measure models}, 
      author={Nihat Ay and Jürgen Jost and Hông Vân Lê and Lorenz Schwachhöfer},
      year={2017},
      eprint={1510.07305},
      archivePrefix={arXiv},
      primaryClass={math.DG},
      url={https://arxiv.org/abs/1510.07305}, 
}

@article{Itoh_2023,
   title={Geometric mean of probability measures and geodesics of Fisher information metric},
   volume={296},
   ISSN={1522-2616},
   url={http://dx.doi.org/10.1002/mana.202000167},
   DOI={10.1002/mana.202000167},
   number={5},
   journal={Mathematische Nachrichten},
   publisher={Wiley},
   author={Itoh, Mitsuhiro and Satoh, Hiroyasu},
   year={2023},
   month=feb, pages={1901–1927} }

@article{Bruveris_2018,
   title={Geometry of the Fisher–Rao metric on the space of smooth densities on a compact manifold},
   volume={292},
   ISSN={1522-2616},
   url={http://dx.doi.org/10.1002/mana.201600523},
   DOI={10.1002/mana.201600523},
   number={3},
   journal={Mathematische Nachrichten},
   publisher={Wiley},
   author={Bruveris, Martins and Michor, Peter W.},
   year={2018},
   month=nov, pages={511–523} }

@misc{cho2025statisticalbergmangeometry,
      title={Statistical Bergman geometry}, 
      author={Gunhee Cho and Jihun Yum},
      year={2025},
      eprint={2305.10207},
      archivePrefix={arXiv},
      primaryClass={math.CV},
      url={https://arxiv.org/abs/2305.10207}, 
}

@article{Facchi_2010,
   title={Classical and quantum Fisher information in the geometrical formulation of quantum mechanics},
   volume={374},
   ISSN={0375-9601},
   url={http://dx.doi.org/10.1016/j.physleta.2010.10.005},
   DOI={10.1016/j.physleta.2010.10.005},
   number={48},
   journal={Physics Letters A},
   publisher={Elsevier BV},
   author={Facchi, Paolo and Kulkarni, Ravi and Man’ko, V.I. and Marmo, Giuseppe and Sudarshan, E.C.G. and Ventriglia, Franco},
   year={2010},
   month=nov, pages={4801–4803} }

@misc{brehmer2020flowssimultaneousmanifoldlearning,
      title={Flows for simultaneous manifold learning and density estimation}, 
      author={Johann Brehmer and Kyle Cranmer},
      year={2020},
      eprint={2003.13913},
      archivePrefix={arXiv},
      primaryClass={stat.ML},
      url={https://arxiv.org/abs/2003.13913}, 
}

@misc{ross2021tractabledensityestimationlearned,
      title={Tractable Density Estimation on Learned Manifolds with Conformal Embedding Flows}, 
      author={Brendan Leigh Ross and Jesse C. Cresswell},
      year={2021},
      eprint={2106.05275},
      archivePrefix={arXiv},
      primaryClass={stat.ML},
      url={https://arxiv.org/abs/2106.05275}, 
}

@misc{kreutzdelgado2009complexgradientoperatorcrcalculus,
      title={The Complex Gradient Operator and the CR-Calculus}, 
      author={Ken Kreutz-Delgado},
      year={2009},
      eprint={0906.4835},
      archivePrefix={arXiv},
      primaryClass={math.OC},
      url={https://arxiv.org/abs/0906.4835}, 
}

@misc{collins2012twistedkahlerricciflow,
      title={The twisted Kahler-Ricci flow}, 
      author={Tristan C. Collins and Gábor Székelyhidi},
      year={2012},
      eprint={1207.5441},
      archivePrefix={arXiv},
      primaryClass={math.DG},
      url={https://arxiv.org/abs/1207.5441}, 
}

@misc{chen2005noteuniformizationriemannsurfaces,
      title={A note on uniformization of Riemann surfaces by Ricci flow}, 
      author={Xiuxiong Chen and Peng Lu and Gang Tian},
      year={2005},
      eprint={math/0505163},
      archivePrefix={arXiv},
      primaryClass={math.DG},
      url={https://arxiv.org/abs/math/0505163}, 
}

@misc{chau2007surveykahlerricciflowyaus,
      title={A survey on the K\"ahler-Ricci flow and Yau's uniformization conjecture}, 
      author={Albert Chau and Luen-Fai Tam},
      year={2007},
      eprint={math/0702257},
      archivePrefix={arXiv},
      primaryClass={math.DG},
      url={https://arxiv.org/abs/math/0702257}, 
}

@misc{chen2002uniformizationtheoremcompletenoncompact,
      title={A Uniformization Theorem Of Complete Noncompact K\"{a}hler Surfaces With Positive Bisectional Curvature}, 
      author={Bing-Long Chen and Siu-Hung Tang and Xi-Ping Zhu},
      year={2002},
      eprint={math/0211372},
      archivePrefix={arXiv},
      primaryClass={math.DG},
      url={https://arxiv.org/abs/math/0211372}, 
}

@misc{chau2003gradientkahlerriccisolitonsuniformization,
      title={Gradient K\"ahler-Ricci solitons and a uniformization conjecture}, 
      author={Albert Chau and Luen-Fai Tam},
      year={2003},
      eprint={math/0310198},
      archivePrefix={arXiv},
      primaryClass={math.DG},
      url={https://arxiv.org/abs/math/0310198}, 
}

@misc{BamlerRicciFlowNotes2015,
  author       = {Bamler, Richard},
  title        = {Ricci Flow Lecture Notes},
  howpublished = {Notes by Otis Chodosh and Christos Mantoulidis},
  year         = {2015},
  month        = {June},
  url          = {https://web.stanford.edu/~ochodosh/Bamler-RFnotes.pdf},
  note         = {Date: June 8, 2015}
}

@book{Topping2006,
  author    = {Topping, Peter},
  title     = {Lectures on the Ricci Flow},
  series    = {London Mathematical Society Lecture Note Series},
  volume    = {325},
  year      = {2006},
  publisher = {Cambridge University Press},
  address   = {Cambridge, UK},
  doi       = {10.1017/CBO9780511721465}
}

@misc{collins2018inversemongeampereflowapplications,
      title={The inverse Monge-Ampere flow and applications to Kahler-Einstein metrics}, 
      author={Tristan C. Collins and Tomoyuki Hisamoto and Ryosuke Takahashi},
      year={2018},
      eprint={1712.01685},
      archivePrefix={arXiv},
      primaryClass={math.DG},
      url={https://arxiv.org/abs/1712.01685}, 
}

@misc{klemyatin2025convergenceinversemongeampereflow,
      title={Convergence of the inverse Monge-Ampere flow and Nadel multiplier ideal sheaves}, 
      author={Nikita Klemyatin},
      year={2025},
      eprint={2411.17978},
      archivePrefix={arXiv},
      primaryClass={math.DG},
      url={https://arxiv.org/abs/2411.17978}, 
}

@misc{perelman2002entropyformularicciflow,
      title={The entropy formula for the Ricci flow and its geometric applications}, 
      author={Grisha Perelman},
      year={2002},
      eprint={math/0211159},
      archivePrefix={arXiv},
      primaryClass={math.DG},
      url={https://arxiv.org/abs/math/0211159}, 
}

@misc{scagliotti2024normalizingflowsapproximationsoptimal,
      title={Normalizing flows as approximations of optimal transport maps via linear-control neural ODEs}, 
      author={Alessandro Scagliotti and Sara Farinelli},
      year={2024},
      eprint={2311.01404},
      archivePrefix={arXiv},
      primaryClass={math.OC},
      url={https://arxiv.org/abs/2311.01404}, 
}

@article{Huang_2023,
   title={Bridging mean-field games and normalizing flows with trajectory regularization},
   volume={487},
   ISSN={0021-9991},
   url={http://dx.doi.org/10.1016/j.jcp.2023.112155},
   DOI={10.1016/j.jcp.2023.112155},
   journal={Journal of Computational Physics},
   publisher={Elsevier BV},
   author={Huang, Han and Yu, Jiajia and Chen, Jie and Lai, Rongjie},
   year={2023},
   month=aug, pages={112155} }

@misc{zhang2023meanfieldgameslaboratorygenerative,
      title={A mean-field games laboratory for generative modeling}, 
      author={Benjamin J. Zhang and Markos A. Katsoulakis},
      year={2023},
      eprint={2304.13534},
      archivePrefix={arXiv},
      primaryClass={stat.ML},
      url={https://arxiv.org/abs/2304.13534}, 
}

@misc{lipman2023flowmatchinggenerativemodeling,
      title={Flow Matching for Generative Modeling}, 
      author={Yaron Lipman and Ricky T. Q. Chen and Heli Ben-Hamu and Maximilian Nickel and Matt Le},
      year={2023},
      eprint={2210.02747},
      archivePrefix={arXiv},
      primaryClass={cs.LG},
      url={https://arxiv.org/abs/2210.02747}, 
}

@misc{papamakarios2021normalizingflowsprobabilisticmodeling,
      title={Normalizing Flows for Probabilistic Modeling and Inference}, 
      author={George Papamakarios and Eric Nalisnick and Danilo Jimenez Rezende and Shakir Mohamed and Balaji Lakshminarayanan},
      year={2021},
      eprint={1912.02762},
      archivePrefix={arXiv},
      primaryClass={stat.ML},
      url={https://arxiv.org/abs/1912.02762}, 
}

@misc{zhai2025normalizingflowscapablegenerative,
      title={Normalizing Flows are Capable Generative Models}, 
      author={Shuangfei Zhai and Ruixiang Zhang and Preetum Nakkiran and David Berthelot and Jiatao Gu and Huangjie Zheng and Tianrong Chen and Miguel Angel Bautista and Navdeep Jaitly and Josh Susskind},
      year={2025},
      eprint={2412.06329},
      archivePrefix={arXiv},
      primaryClass={cs.CV},
      url={https://arxiv.org/abs/2412.06329}, 
}

@misc{grathwohl2018ffjordfreeformcontinuousdynamics,
      title={FFJORD: Free-form Continuous Dynamics for Scalable Reversible Generative Models}, 
      author={Will Grathwohl and Ricky T. Q. Chen and Jesse Bettencourt and Ilya Sutskever and David Duvenaud},
      year={2018},
      eprint={1810.01367},
      archivePrefix={arXiv},
      primaryClass={cs.LG},
      url={https://arxiv.org/abs/1810.01367}, 
}

@inbook{inbook,
author = {Villani, Cédric},
year = {2008},
month = {01},
pages = {xxii+973},
title = {Optimal transport -- Old and new},
volume = {338},
doi = {10.1007/978-3-540-71050-9}
}

@misc{kingma2018glowgenerativeflowinvertible,
      title={Glow: Generative Flow with Invertible 1x1 Convolutions}, 
      author={Diederik P. Kingma and Prafulla Dhariwal},
      year={2018},
      eprint={1807.03039},
      archivePrefix={arXiv},
      primaryClass={stat.ML},
      url={https://arxiv.org/abs/1807.03039}, 
}

@misc{dinh2017densityestimationusingreal,
      title={Density estimation using Real NVP}, 
      author={Laurent Dinh and Jascha Sohl-Dickstein and Samy Bengio},
      year={2017},
      eprint={1605.08803},
      archivePrefix={arXiv},
      primaryClass={cs.LG},
      url={https://arxiv.org/abs/1605.08803}, 
}

@misc{tran2019discreteflowsinvertiblegenerative,
      title={Discrete Flows: Invertible Generative Models of Discrete Data}, 
      author={Dustin Tran and Keyon Vafa and Kumar Krishna Agrawal and Laurent Dinh and Ben Poole},
      year={2019},
      eprint={1905.10347},
      archivePrefix={arXiv},
      primaryClass={cs.LG},
      url={https://arxiv.org/abs/1905.10347}, 
}

@inproceedings{NEURIPS2022_e9e1a0ab,
 author = {Huang, Hengguan and Gu, Xiangming and Wang, Hao and Xiao, Chang and Liu, Hongfu and Wang, Ye},
 booktitle = {Advances in Neural Information Processing Systems},
 editor = {S. Koyejo and S. Mohamed and A. Agarwal and D. Belgrave and K. Cho and A. Oh},
 pages = {36000--36013},
 publisher = {Curran Associates, Inc.},
 title = {Extrapolative Continuous-time Bayesian Neural Network for Fast Training-free Test-time Adaptation},
 url = {https://proceedings.neurips.cc/paper_files/paper/2022/file/e9e1a0abc1a5b19a4aeb80dab19c82ae-Paper-Conference.pdf},
 volume = {35},
 year = {2022}
}

@inproceedings{NEURIPS2020_df1a336b,
 author = {Suzuki, Taiji},
 booktitle = {Advances in Neural Information Processing Systems},
 editor = {H. Larochelle and M. Ranzato and R. Hadsell and M.F. Balcan and H. Lin},
 pages = {19224--19237},
 publisher = {Curran Associates, Inc.},
 title = {Generalization bound of globally optimal non-convex neural network training: Transportation map estimation by infinite dimensional Langevin dynamics},
 url = {https://proceedings.neurips.cc/paper_files/paper/2020/file/df1a336b7e0b0cb186de6e66800c43a9-Paper.pdf},
 volume = {33},
 year = {2020}
}

@misc{zhang2021diffusionnormalizingflow,
      title={Diffusion Normalizing Flow}, 
      author={Qinsheng Zhang and Yongxin Chen},
      year={2021},
      eprint={2110.07579},
      archivePrefix={arXiv},
      primaryClass={cs.LG},
      url={https://arxiv.org/abs/2110.07579}, 
}

@misc{shen2022cherncalabiflowhermitianmanifolds,
      title={A Chern-Calabi flow on Hermitian manifolds}, 
      author={Xi Sisi Shen},
      year={2022},
      eprint={2011.09683},
      archivePrefix={arXiv},
      primaryClass={math.DG},
      url={https://arxiv.org/abs/2011.09683}, 
}

@misc{perelman2003ricciflowsurgerythreemanifolds,
      title={Ricci flow with surgery on three-manifolds}, 
      author={Grisha Perelman},
      year={2003},
      eprint={math/0303109},
      archivePrefix={arXiv},
      primaryClass={math.DG},
      url={https://arxiv.org/abs/math/0303109}, 
}

@misc{calabi2001spacekahlermetricsii,
      title={The Space of K\"ahler metrics (II)}, 
      author={E. Calabi and X. X. Chen},
      year={2001},
      eprint={math/0108162},
      archivePrefix={arXiv},
      primaryClass={math.DG},
      url={https://arxiv.org/abs/math/0108162}, 
}

@misc{he2025twistedcalabifunctionaltwisted,
      title={Twisted Calabi functional and twisted Calabi flow}, 
      author={Jie He and Haozhao Li},
      year={2025},
      eprint={2512.02451},
      archivePrefix={arXiv},
      primaryClass={math.DG},
      url={https://arxiv.org/abs/2512.02451}, 
}

@article{Cao1985,
  author    = {Cao, Huai-Dong},
  title     = {Deformation of {K}ähler metrics to {K}ähler--{E}instein metrics on compact {K}ähler manifolds},
  journal   = {Inventiones mathematicae},
  year      = {1985},
  volume    = {81},
  number    = {2},
  pages     = {359--372},
  doi       = {10.1007/BF01389058},
  url       = {https://doi.org/10.1007/BF01389058},
  issn      = {1432-1297},
}

@misc{phong2007multiplieridealsheaveskahlerricci,
      title={Multiplier Ideal Sheaves and the K\"ahler-Ricci Flow}, 
      author={D. H. Phong and Natasa Sesum and Jacob Sturm},
      year={2007},
      eprint={math/0611794},
      archivePrefix={arXiv},
      primaryClass={math.DG},
      url={https://arxiv.org/abs/math/0611794}, 
}

@misc{chen2009stabilitykahlerricciflow,
      title={Stability of K\"ahler-Ricci flow}, 
      author={Xiuxiong Chen and Haozhao Li},
      year={2009},
      eprint={0801.3086},
      archivePrefix={arXiv},
      primaryClass={math.DG},
      url={https://arxiv.org/abs/0801.3086}, 
}

@article{Huang_2020,
   title={Kähler–Ricci flow on homogeneous toric bundles},
   volume={31},
   ISSN={1793-6519},
   url={http://dx.doi.org/10.1142/S0129167X20500226},
   DOI={10.1142/s0129167x20500226},
   number={03},
   journal={International Journal of Mathematics},
   publisher={World Scientific Pub Co Pte Ltd},
   author={Huang, Hong},
   year={2020},
   month=feb, pages={2050022} }

@misc{salman2018deepdiffeomorphicnormalizingflows,
      title={Deep Diffeomorphic Normalizing Flows}, 
      author={Hadi Salman and Payman Yadollahpour and Tom Fletcher and Kayhan Batmanghelich},
      year={2018},
      eprint={1810.03256},
      archivePrefix={arXiv},
      primaryClass={stat.ML},
      url={https://arxiv.org/abs/1810.03256}
}

@article{Calamai_2015,
   title={The Dirichlet and the weighted metrics for the space of Kähler metrics},
   volume={363},
   ISSN={1432-1807},
   url={http://dx.doi.org/10.1007/s00208-015-1188-x},
   DOI={10.1007/s00208-015-1188-x},
   number={3–4},
   journal={Mathematische Annalen},
   publisher={Springer Science and Business Media LLC},
   author={Calamai, Simone and Zheng, Kai},
   year={2015},
   month=mar, pages={817–856} }

@article{Chen_2013,
   title={The pseudo-Calabi flow},
   volume={2013},
   ISSN={0075-4102},
   url={http://dx.doi.org/10.1515/crelle.2012.033},
   DOI={10.1515/crelle.2012.033},
   number={674},
   journal={Journal für die reine und angewandte Mathematik (Crelles Journal)},
   publisher={Walter de Gruyter GmbH},
   author={Chen, Xiuxiong and Zheng, Kai},
   year={2013},
   month=jan }

@article{Roch_2026,
   title={Learning informed prior distributions with normalizing flows for Bayesian analysis},
   volume={113},
   ISSN={2469-9993},
   url={http://dx.doi.org/10.1103/wms4-m4cb},
   DOI={10.1103/wms4-m4cb},
   number={3},
   journal={Physical Review C},
   publisher={American Physical Society (APS)},
   author={Roch, Hendrik and Shen, Chun},
   year={2026},
   month=mar }

@misc{rosenhahn2024quantumnormalizingflowsanomaly,
      title={Quantum Normalizing Flows for Anomaly Detection}, 
      author={Bodo Rosenhahn and Christoph Hirche},
      year={2024},
      eprint={2402.02866},
      archivePrefix={arXiv},
      primaryClass={quant-ph},
      url={https://arxiv.org/abs/2402.02866}, 
}

@misc{shoham2025flatnessall,
      title={Flatness After All?}, 
      author={Neta Shoham and Liron Mor-Yosef and Haim Avron},
      year={2025},
      eprint={2506.17809},
      archivePrefix={arXiv},
      primaryClass={cs.LG},
      url={https://arxiv.org/abs/2506.17809}, 
}
\endgroup

\appendix

\section{Notation}

\begin{table}[H]
\centering
\renewcommand{\arraystretch}{1.0} 
\small
\begin{tabular}{@{} l p{11cm} @{}}
\toprule
\textbf{Symbol} & \textbf{Description} \\ \midrule
\rowcolor{NavajoWhite!15} $\Phi$ & Kähler potential \\
\addlinespace[6pt] \rowcolor{Maroon!15} $\dot{\Phi}$ & Time derivative of the Kähler potential \\
\addlinespace[6pt] \rowcolor{NavajoWhite!15} $h$ & Kähler (Hermitian) metric \\
\addlinespace[6pt] \rowcolor{Maroon!15} $\text{Ric}_{i \overline{j}}$ & Complex Ricci curvature \\
\addlinespace[6pt] \rowcolor{NavajoWhite!15} $R$ & Scalar curvature \\
\addlinespace[6pt] \rowcolor{Maroon!15} $\partial_i \partial_{\overline{j}}$ & Mixed second-order Wirtinger derivatives with respect to $z^i, \overline{z}^j$ \\
\addlinespace[6pt] \rowcolor{NavajoWhite!15} $\text{div}_{\omega_t}$ & Divergence as in section \ref{sec:notations} \\
\addlinespace[6pt] \rowcolor{Maroon!15} $\nabla, \nabla_z$ & Vectorized Wirtinger derivative / Wirtinger Jacobian \\
\addlinespace[6pt] \rowcolor{NavajoWhite!15} $\Psi$ & Flow map of the complex normalizing flow \\
\addlinespace[6pt] \rowcolor{Maroon!15} $\text{KL}(q_t \parallel p)$ & Kullback-Leibler divergence $\int \log \left( \frac{dq_t}{dp} \right) dq_t$ \\ 
\addlinespace[6pt] \rowcolor{NavajoWhite!15} $(\partial, \overline{\partial})$ & Dolbeault operators \\
\addlinespace[6pt] \rowcolor{Maroon!15} $\Delta_{\overline{\partial}} f$ & Kähler Laplacian operator $h^{\overline{j} i} \partial_{i} \partial_{\overline{j}} f$\\
\addlinespace[6pt] \rowcolor{NavajoWhite!15} $\omega_t^d$ & Kähler top form of $\sqrt{-1} \partial \overline{\partial} \Phi$ at time $t$ \\
\addlinespace[6pt] \rowcolor{Maroon!15} $\langle \alpha, \beta \rangle_{\omega_t}$ & Hermitian inner product, in local coordinates $h^{\overline{j} i} \alpha_i \beta_{\overline{j}}$ \\
\addlinespace[6pt] \rowcolor{NavajoWhite!15} $\langle \partial f, \partial g \rangle_{\omega_t}$ & Hermitian inner product of differentials  $h^{\overline{j} i} \partial_i f \partial_{\overline{j}} g$  \\
\addlinespace[6pt] \rowcolor{Maroon!15}  $| \partial \overline{\partial} f|^2$ & $|\partial \overline{\partial} f|^2 = h^{\overline{\ell} i} h^{\overline{j} k} (\partial_i \partial_{\overline{j}} f) (\partial_k \partial_{\overline{\ell}} f)$ \\
\addlinespace[6pt] \rowcolor{NavajoWhite!15}  $| \partial f|^2$ & $h^{\overline{j} i } \partial_i f \partial_{\overline{j}} f$ \\
\addlinespace[6pt] \rowcolor{Maroon!15} $ \nabla_i, \nabla^i$ & Levi-Civita connection and its index-raised version
\\
\addlinespace[6pt] \rowcolor{NavajoWhite!15}$\text{Tr}_{\omega_t}(\alpha)$ & Trace of the (1,1)-form $\alpha$ $h^{\overline{j} i} \alpha_{i \overline{j}}$ \\
\addlinespace[6pt] \rowcolor{Maroon!15}   $\doublenabla$ & Gradient w.r.t. the Dirichlet metric  \\
\addlinespace[6pt] \rowcolor{NavajoWhite!15} $ \psi$ & Twisting potential \\
\bottomrule
\end{tabular}
\end{table}

\section{Vanishing Ricci curvature under holomorphic pullbacks}
\label{app:vanishing_ricci}

For a holomorphic map $\Psi: M \to \mathbb{C}^n$ between complex manifolds of the same dimension, the metric components are given by the metric $h_{i\overline{j}} = \sum_{k} \frac{\partial \Psi^k}{\partial z^i} \overline{\frac{\partial \Psi^k}{\partial z^j}}$. In matrix notation, $h = (\nabla \Psi)^{\dagger} (\nabla \Psi)$, where $\nabla \Psi$ is the holomorphic (Wirtinger) Jacobian. The determinant of the metric tensor is
\begin{align}
\det(h) = \det(\nabla \Psi) \det(\overline{\nabla \Psi}) = |\det(\nabla \Psi)|^2 .
\end{align}
The Ricci form $\rho$ for a Kähler metric is defined as
\begin{align}
\rho = -\sqrt{-1} \partial \overline{\partial} \log \det(h) .
\end{align}
Substituting the expression for the determinant
\begin{align}
\rho = -\sqrt{-1} \partial \overline{\partial} \log \left( \det(\nabla \Psi) \overline{\det(\nabla \Psi)} \right) = -\sqrt{-1} \partial \overline{\partial} \left( \log \det(\nabla \Psi) + \log \overline{\det(\nabla \Psi)} \right) .
\end{align}
Since $\Psi$ is holomorphic, consequently, $\overline{\partial} \log J = 0, \partial \log \overline{J} = 0$. In particular, one term vanishes under $\partial$ and the other vanishes under $\overline{\partial}$.

\section{Change of variables identity}
\label{app:change_of_variables}

For a complex manifold of dimension $d$ with Hermitian metric $h_{i\overline{j}}$, the volume form $\Omega$ obeys
\begin{align}
\Omega = \det(h) \left( \frac{\sqrt{-1}}{2} \right)^d dz^1 \wedge d\overline{z}^1 \wedge \hdots \wedge dz^d \wedge d\overline{z}^d .
\end{align}
Let us specialize to the normalizing flow. We have a sequence of measures $q_k$ associated with metrics $h_k$. If we assume the density $q$ is identified with the metric volume density, we have
\begin{align}
d\mu_k = \det(h_k) dV
\end{align}
where $dV = (\frac{\sqrt{-1}}{2})^d \bigwedge_j dz^j \wedge d\overline{z}^j$ is the standard Euclidean volume element. Consider the map $\Psi_k: M_{k-1} \to M_k$ where $z_k = \Psi_k(z_{k-1})$. Even though $\Psi_k$ is non-holomorphic, it is a smooth diffeomorphism. The transformation of the coordinate volume element $dV$ under $\Psi_k$ is governed by the real Jacobian. For the complex augmented Jacobian $\mathcal{J}_k$,
\begin{align}
\mathcal{J}_k = \begin{pmatrix} \partial \Psi & \overline{\partial} \Psi \\ \partial \overline{\Psi} & \overline{\partial} \overline{\Psi} \end{pmatrix} .
\end{align}
The relation to the real Jacobian $J_{\mathbb{R}}$ is $\det J_{\mathbb{R}} = |\det \mathcal{J}_k|$. Consequently, the volume element transforms as
\begin{align}
\Psi_k^*(dV_k) = |\det \mathcal{J}_k| dV_{k-1} .
\end{align}
The fundamental law of normalizing flows (i.e., conservation of probability mass) states that for any measurable set
\begin{align}
\int_A q_{k-1}(z_{k-1}) dV_{k-1} = \int_{\Psi_k(A)} q_k(z_k) dV_k .
\end{align}
By the change of variables formula, the right-hand side becomes
\begin{align}
\int_A q_k(\Psi_k(z_{k-1})) |\det \mathcal{J}_k| dV_{k-1} .
\end{align}
Since $A$ measurable is arbitrary, we may equate the integrands in a fundamental lemma of calculus of variations-type manner
\begin{align}
q_{k-1}(z_{k-1}) = q_k(z_k) |\det \mathcal{J}_k| .
\end{align}
We substitute $q \propto \det h$ (which follows by our defining geometry as in section \ref{sec:normalizing_flow_geometry})
\begin{align}
\det h_{k-1} = \det h_k |\det \mathcal{J}_k| .
\end{align}

\section{Connections to normalized equations}
\label{app:normalized_kahler_ricci_flow}

\noindent \textbf{Theorem 2 (continued).} \textit{ The density evolution under our condition but normalized obeys \begin{align}
\log q_{k,\theta} = \log q_{k-1,\theta} - \Big[ \log |\det(\mathcal{J}_{k,\theta})| + \nu \Delta t \log q_{k-1,\theta} \Big]  .
\end{align}
}

\vspace{2mm}

\noindent Our work in this section has partial connection to the normalized Kähler-Ricci flow which evolves, as a form, according to
\begin{align}
\frac{\partial \omega}{\partial t} = -\text{Ric}(\omega) - \nu \omega ,
\end{align}
or as a cohomology class evolution $\frac{\partial}{\partial t} [\omega_t] = -2\pi c_1(M) - \nu [\omega_t]$, where $c_1(M)$ is the first Chern class. In local coordinates, with form $\omega = \sqrt{-1} \sum_{jk} h_{j\overline{k}} dz^j \wedge d\overline{z}^k$, we have 
\begin{align}
\partial_t h_{i \overline{j}} = - \text{Ric}_{i \overline{j}} - \nu h_{i \overline{j}} .
\end{align}
From a Fisher-Bayesian information-theoretic formulation in the continuum limit, we examine
\begin{align}
\lim_{\Delta t \rightarrow 0^+} \frac{I_{i \overline{j}, k} - I_{i \overline{j}, k-1}}{\Delta t} =  -\EX_{\theta \sim p(\theta|\mathcal{D})} \Big[ \partial_t \text{Ric}_{ i \overline{j}} (h_{k-1,\theta}) + \nu I_{i \overline{j}, k-1} \Big] ,
\end{align}
and reversing the chain as in section \ref{sec:main_contribution}, this corresponds to the normalizing flow, preceding the mixed Wirtinger differentiation,
\begin{align}
\log q_{K,\theta}(z) = \log q_{0,\theta}(z) - \sum_{k=1}^K \Big[ \log | \det ( \mathcal{J}_{k,\theta}) | + \nu \Delta t \log q_{k,\theta}(z) \Big] .
\end{align}
We note that we need the same signs in the last two terms. When we take Ricci curvature, a sign is included. Also, the Fisher information defined via second order also includes a sign change. An additional negative is added upon switching the Fisher information order. We motivate the following approximation by the ansatz $\Gamma$ (if it exists) such that the full Jacobian satisfies
\begin{align}
\mathcal{J}_{k,\theta}^\Gamma = [q_{k-1,\theta}(\Psi_{k,\theta})]^{\frac{\nu \Delta t}{2d}} \mathcal{J}_{k,\theta}^\Psi + O(\Delta t) .
\end{align}
Taking the absolute value of the determinant and factoring the scalar out of the $2d \times 2d$ matrix gives
\begin{align}
|\det(\mathcal{J}_{k,\theta}^\Gamma)| & \approx \Big| \det \Big( [q_{k-1,\theta}(\Psi_{k,\theta})]^{\frac{\nu \Delta t}{2d}} \mathcal{J}_{k,\theta}^\Psi \Big) \Big| \\
& = [q_{k-1,\theta}(\Psi_{k,\theta})]^{\frac{\nu \Delta t}{2d} \cdot 2d} |\det(\mathcal{J}_{k,\theta}^\Psi)| \\
& = q_{k-1,\theta}(\Psi_{k,\theta})^{\nu \Delta t} |\det(\mathcal{J}_{k,\theta}^\Psi)| .
\end{align}
Taking the log,
\begin{align}
\log |\det(\mathcal{J}_{k,\theta}^\Gamma)| \approx \log |\det(\mathcal{J}_{k,\theta}^\Psi)| + \nu \Delta t \log q_{k-1,\theta}(\Psi_{k,\theta}) .
\end{align}
By change of variables,
\begin{align}
\log q_{k,\theta}(\Gamma_{k}) = \log q_{k-1,\theta}(\Gamma_{k-1}) - \log|\det(\mathcal{J}_{k,\theta}^\Gamma)| .
\end{align}
Substituting gives
\begin{align}
\log q_{k,\theta}(z) = \log q_{k-1,\theta}(z) - \Big[ \log |\det \mathcal{J}_{k,\theta}(z) | + \nu \Delta t \log q_{k-1,\theta}(z) \Big] 
\end{align}
as desired.

\section{Statistical density of the volume element and Ricci curvature evolution}
\label{app:statistical_densities}

In these sections, we provide statistical evolution quantities relating to classical geometric features, i.e. Ricci curvature and scalar curvature, and their equivalent formulations under the complex normalizing flow. The complex normalizing flow is primarily a statistical-driven archetype, innate to machine learning, thus traditional geometric qualities can be related under their baseline geometric definitions to their statistical varieties under the abridged complex normalizing flow to its Kähler-Ricci variants.

\vspace{2mm}

\noindent \textbf{Theorem 1 (continued).} \textit{Assume local coordinates are time-independent. Then Ricci curvature statistically obeys
\begin{align}
\text{Ric}_{i \overline{j}} = - \EX_p \Big[ \text{div}_{\mathbb{R}^n}(f)(\partial_i \partial_{\overline{j}} \log q ) \Big] + \EX_p \Big[ \partial_i \partial_{\overline{j}} ( \text{div}_{\omega_t}(g) + \partial_t \log \omega_t) \Big] ,
\end{align}
and moreover, the Ricci curvature time derivative obeys
\begin{align}
\frac{\partial}{\partial t} \text{Ric}_{\ell \overline{k}} = -\partial_\ell \partial_{\overline{k}} \Bigg( h^{\overline{j} i} \Big( \EX_p \Big[ \text{div}_{\mathbb{R}^n}(f)(\partial_i \partial_{\overline{j}} \log q ) \Big] + \EX_p \Big[ \partial_i \partial_{\overline{j}} ( \text{div}_{\omega_t}(g) + \partial_t \log \omega_t) \Big] \Big) \Bigg) .
\end{align}
}

\vspace{2mm}

\noindent \textit{Proof.} Again denoting $p$ the Bayesian density of the parameter and $q$ the density of the empirical data, we derive the following results. Let us first note the complex instantaneous change of variables identities \cite{chen2019neuralordinarydifferentialequations}
\begin{align}
\begin{cases}
& \frac{d}{d t} \log p = - \text{div}_{\mathbb{R}^n}(f) 
\\
& \frac{d}{d t} \log q = - \text{div}_{\omega_t}(g) - \partial_t \log \omega_t .
\end{cases}
\end{align}
Let us revisit our information metric
\begin{align}
h_{i \overline{j}} = \EX_p \Big[ -\partial_i \partial_{\overline{j}} \log q \Big] .
\end{align}
Let us reformulate this an integral and differentiating,
\begin{align}
\frac{\partial}{\partial t} h_{i \overline{j}} = \frac{\partial }{\partial t} \int_{\Theta} p(\theta,t)( - \partial_i \partial_{\overline{j}} \log q_{\theta}(z,t) ) d\theta .
\end{align}
Using shorthand notation for $p(\theta,t), q_{\theta}(z,t)$, via product rule,
\begin{align}
\frac{\partial}{\partial t} h_{i \overline{j}} & = \int_{\Theta} \dot{p}(- \partial_i \partial_{\overline{j}} \log q )d\theta + \int_{\Theta} p (- \partial_i \partial_{\overline{j}} \frac{\partial}{\partial t} \log q )d\theta
\\
& = \int_{\Theta} -p\text{div}_{\mathbb{R}^n}(f) ( - \partial_i \partial_{\overline{j}} \log q) d\theta + \int_{\Theta} p (-\partial_i \partial_{\overline{j}}(- \text{div}_{\omega_t}(g) - \partial_t \log \omega_t)) d\theta
\\ 
& = \EX_p \Big[ \text{div}_{\mathbb{R}^n}(f)(\partial_i \partial_{\overline{j}} \log q ) \Big] + \EX_p \Big[ \partial_i \partial_{\overline{j}} ( \text{div}_{\omega_t}(g) + \partial_t \log \omega_t)  \Big] .
\end{align}
It is beneficial to keep in mind $\theta$ is not complex-valued, so the above is ordinary Lebesgue measure. By Jacobi's formula,
\begin{align}
\frac{\partial}{\partial t} \log \det ( h ) = \text{Tr}(h^{-1} \dot{h}) = h^{\overline{j} i} \dot{h}_{i \overline{j}} .
\end{align}
Thus,
\begin{align}
\frac{\partial}{\partial t} \log \det ( h ) = h^{\overline{j} i} \Big( \EX_p \Big[ \text{div}_{\mathbb{R}^n}(f)(\partial_i \partial_{\overline{j}} \log q ) \Big] + \EX_p \Big[ \partial_i \partial_{\overline{j}} ( \text{div}_{\omega_t}(g) + \partial_t \log \omega_t) \Big] \Big) .
\end{align}
Therefore, we can compute the time derivative of Ricci curvature as
\begin{align}
\frac{\partial}{\partial t} \text{Ric}_{\ell \overline{k}} = -\partial_{\ell} \partial_{\overline{k}} \Bigg( h^{\overline{j} i} \Big( \EX_p \Big[ \text{div}_{\mathbb{R}^n}(f)(\partial_i \partial_{\overline{j}} \log q ) \Big] + \EX_p \Big[ \partial_i \partial_{\overline{j}} ( \text{div}_{\omega_t}(g) + \partial_t \log \omega_t) \Big] \Big) \Bigg) .
\end{align}
Notice the above can be simplified greatly using the Kähler-Ricci flow, therefore we have found a bridge between the instantaneous change of variables formula and Kähler-Ricci flow. Using the flow identity $\partial_t h_{i \overline{j}} = -\text{Ric}_{i \overline{j}}$, Jacobi's formula becomes
\begin{align}
\label{eqn:jacobis_scalar}
\frac{\partial}{\partial t} \log \det ( h ) = \text{Tr}(h^{-1} \dot{h}) = h^{\overline{j} i} \dot{h}_{i \overline{j}} = -h^{\overline{j} i} \EX \partial_t \text{Ric}_{i \overline{j}}.
\end{align}
Moreover, 
\begin{align}
\text{Ric}_{i \overline{j}} = - \EX_p \Big[ \text{div}_{\mathbb{R}^n}(f)(\partial_i \partial_{\overline{j}} \log q ) \Big] + \EX_p \Big[ \partial_i \partial_{\overline{j}} ( \text{div}_{\omega_t}(g) + \partial_t \log \omega_t) \Big] .
\end{align}

\subsection{Statistical evolution of scalar curvature}

\noindent \textbf{Theorem 1 (continued).}  \textit{ The time derivative of scalar curvature obeys
\begin{align}
\frac{\partial}{\partial t} R  & = - h^{\overline{j} \ell} h^{\overline{k} i} \partial_t \text{Ric}_{\ell \overline{k}} \text{Ric}_{i \overline{j}} + h^{\overline{j} i} \partial_i \partial_{\overline{j}} ( h^{\overline{\ell} k} \text{Ric}_{k \overline{\ell}} ) ,
\end{align}
which is consistent with \cite{song2012lecturenoteskahlerricciflow}.
}

\vspace{2mm}

\noindent \textit{Proof.} By definition $R = h^{\overline{j} i} \text{Ric}_{i \overline{j}}$. Differentiating scalar curvature,
\begin{align}
\frac{\partial R}{\partial t} & = \frac{\partial h^{\overline{j} i}}{\partial t} \text{Ric}_{i \overline{j}} + h^{\overline{j} i} \frac{\partial \text{Ric}_{i \overline{j}}}{\partial t} 
\\
&  = -(h^{\overline{j} \ell} h^{\overline{k} i} \partial_t \EX \text{Ric}_{\ell \overline{k}}) \text{Ric}_{i \overline{j}} - h^{\overline{j} i} \partial_i \partial_{\overline{j}} \Bigg( h^{\overline{\ell} k} \Big( \EX_p \Big[ \text{div}_{\mathbb{R}^n}(f)(\partial_k \partial_{\overline{\ell}} \log q ) \Big] + \EX_p \Big[ \partial_k \partial_{\overline{\ell}} \text{div}_{\omega}(g) \Big] \Big) \Bigg)
\\
& = -h^{\overline{j} \ell} h^{\overline{k} i } \Bigg( -\partial_t \EX_p \Big[ \text{div}_{\mathbb{R}^n}(f)(\partial_{\ell} \partial_{\overline{k}} \log q ) \Big] + \partial_t \EX_p \Big[ \partial_\ell \partial_{\overline{k}} ( \text{div}_{\omega_t}(g) + \partial_t \log \omega_t) \Big]  \Bigg) 
\\
& \ \ \ \ \ \ \ \ \ \ \ \ \times \Bigg( - \EX_p \Big[ \text{div}_{\mathbb{R}^n}(f)(\partial_i \partial_{\overline{j}} \log q ) \Big] + \EX_p \Big[ \partial_i \partial_{\overline{j}} ( \text{div}_{\omega_t}(g) + \partial_t \log \omega_t) \Big] \Bigg) 
\\
& \ \ \ \ \ \ \ \ \ \ \ \ -  h^{\overline{j} i} \partial_i \partial_{\overline{j}} \Bigg( h^{\overline{\ell} k} \Big( \EX_p \Big[ \text{div}_{\mathbb{R}^n}(f)(\partial_k \partial_{\overline{\ell}} \log q ) \Big] + \EX_p \Big[ \partial_k \partial_{\overline{\ell}} ( \text{div}_{\omega_t}(g) + \partial_t \log \omega_t) \Big] \Big) \Bigg) .
\end{align}
It can be noted
\begin{align}
\underbrace{\EX_p \Big[ \text{div}_{\mathbb{R}^n}(f)(\partial_i \partial_{\overline{j}} \log q ) \Big]}_{\text{advection}} + \underbrace{\EX_p \Big[ \partial_i \partial_{\overline{j}} ( \text{div}_{\omega_t}(g) + \partial_t \log \omega_t) \Big]}_{\text{diffusion}} .
\end{align}
Now, we get
\begin{align}
\frac{\partial}{\partial t} R  & = -h^{\overline{j} \ell} h^{\overline{k} i} \partial_t \text{Ric}_{\ell \overline{k}} \text{Ric}_{i \overline{j}} - h^{\overline{j} i} \partial_i \partial_{\overline{j}} ( h^{\overline{\ell} k} \text{Ric}_{k \overline{\ell}} )  ,
\end{align}
which is consistent with \cite{song2012lecturenoteskahlerricciflow}.

\section{Relations to optimal transport}
\label{app:continuity_equations}

In this section, we develop connections between classical OT problems with soft and hard-weighted KL divergence weights and kinetic-energy type functional reformulations via the Kähler potential. Let us consider the real-valued optimal transport problem
\begin{align}
\begin{cases}
\min_{\theta} \int_0^T \int_{\Omega} |v(\mathcal{X}(x,t), t; \theta) |^2 \rho_0(x) dx dt
+ \lambda \text{KL}( \rho_T \parallel \mathcal{X}(\cdot,T) \# \rho_0 )
\\
\frac{d \mathcal{X}(x,t)}{dt} = v(\mathcal{X}, t; \theta)
\\
\mathcal{X}(x,0) = x\sim \rho_0 ,
\end{cases}
\end{align}
where $\rho \in \mathcal{P}(\mathbb{R}^d)$,  $\Omega \subseteq \mathbb{R}^d, v: \mathbb{R}^d \times \mathbb{R}^+ \times \Theta \rightarrow \mathbb{R}^d, \mathcal{X} : \mathbb{R}^d \times \mathbb{R}^+ \rightarrow \mathbb{R}^d$. This is a soft KL-weighted constraint problem. Here $\mathcal{X}$ is the Lagrangian flow map and $\rho_0$ is an initial measure in which the starting data derives. A complex normalizing flow setting is
\begin{align}
\begin{cases}
\min_{\theta} \int_0^T \int_{M} |v(\mathcal{Z}(z,t), t; \theta)|^2_{\omega_t} \rho_0(z) \det h \left( \frac{\sqrt{-1}}{2} \right)^d dz_1 \wedge d\overline{z}_1 \wedge \dots \wedge dz_d \wedge d\overline{z}_d dt + \lambda \text{KL}(\rho_T \parallel \mathcal{Z}(\cdot, T) \# \rho_0) \\
\frac{d\mathcal{Z}}{dt} = v(\mathcal{Z}, t; \theta) \\
\mathcal{Z}(z,0) = z\sim \rho_0 \in \mathcal{P}^2(\mathbb{C}^d) .
\end{cases}
\end{align}
The vector field $v$ is a complex vector field of a complex potential, leading to a complex Monge-Ampère equation \cite{zhang2018mongeampereflowgenerativemodeling}. With hard constraints, this is
\begin{align}
\mathcal{L}(\theta, \lambda) = & \int_0^T \int_{M} |v(\mathcal{Z}(z,t), t; \theta)|^2_{\omega_t} \rho_0(z) \left( \frac{\sqrt{-1}}{2} \right)^d dz_1 \wedge d\overline{z}_1 \wedge \dots \wedge dz_d \wedge d\overline{z}_d  dt 
\\
& + \langle \lambda, \rho_T - \mathcal{Z}(\cdot, T) \# \rho_0 \rangle_{C_0(M) \times \mathcal{M}(M)}
\end{align}
subject to optimization
\begin{align}
\begin{cases}
\min_{\theta} \max_{\lambda} \mathcal{L}(\theta, \lambda)
\\
\dfrac{d\mathcal{Z}}{dt} = v(\mathcal{Z}, t; \theta) \\
\mathcal{Z}(z, 0) = z \sim \rho_0 \in \mathcal{P}^2(M) \\
\rho_T = \mathcal{Z}(\cdot, T){\#} \rho_0 .
\end{cases} 
\end{align}
We are using notation
\begin{align}
\langle \lambda, \rho_T - \mathcal{Z}(\cdot, T) \# \rho_0 \rangle_{C_0(M) \times \mathcal{M}(M)} =   \int_{M} \lambda d( \rho_T - \mathcal{Z}(\cdot, T) \# \rho_0), \ \ \ \ \ \lambda \in C_0(M), \rho_T, \rho_0 \in \mathcal{P}^2(M)  , 
\end{align}
and 
\begin{align}
\mathcal{P}^2(M) = \Bigg\{ \mu : \mathcal{B}(M) \rightarrow \mathbb{R}^+ : \int_{M} \mu(dz) = 1 \ \text{and} \ \int_{M} z^{\dagger} z \mu(dz) < \infty \Bigg\} .
\end{align}
It can be noted $\rho_T - \mathcal{Z}(\cdot, T) \# \rho_0$ is a signed measure, and assume that $\mathcal{Z}$ is Borel measurable (we will typically allow it to be smooth). $\mathcal{B}(M)$ is the Borel $\sigma$-algebra of the complex manifold. Elaborating more rigorously in Appendix \ref{app:kinetic_energies}, the above has an equivalent complex potential formulation
\begin{align}
\begin{cases}
\min_{\Phi} \int_0^T \int_{M}  h^{\overline{j}i} \partial_i \dot{\Phi} \partial_{\overline{j}} \dot{\Phi}  \rho_0(z)  \left( \frac{\sqrt{-1}}{2} \right)^d dz_1 \wedge d\overline{z}_1 \wedge \dots \wedge dz_d \wedge d\overline{z}_d dt + \lambda \text{KL}( \rho_T \parallel \mathcal{Z}(\cdot,T) \# \rho_0 )
\\
h_{i\overline{j}} = \partial_i \partial_{\overline{j}} \Phi
\\
\frac{d \mathcal{Z}^i}{dt} = h^{\overline{j} i} \partial_{\overline{j}} \dot{\Phi} .
\end{cases}
\end{align}
The above time derivative is the musical isomorphism. With hard constraints, we get
\begin{align}
\mathcal{L}(\Phi, \lambda) = & \int_0^T \int_{M} h^{\overline{j}i} \partial_i \dot{\Phi} \partial_{\overline{j}} \dot{\Phi} \rho_0(z)  \left( \frac{\sqrt{-1}}{2} \right)^d dz_1 \wedge d\overline{z}_1 \wedge \dots \wedge dz_d \wedge d\overline{z}_d  dt 
\\
& + \langle \lambda, \rho_T - \mathcal{Z}(\cdot, T) \# \rho_0 \rangle_{C_0(M) \times \mathcal{M}(M)}     ,
\end{align} 
with governing dynamics
\begin{align}
\begin{cases}
\min_{\Phi} \max_{\lambda} \mathcal{L}(\Phi, \lambda)
\\
h_{i\overline{j}} = \partial_i \partial_{\overline{j}} \Phi \\
\dfrac{d\mathcal{Z}^i}{dt} = h^{\overline{j}i} \partial_{\overline{j}} \dot{\Phi} \\
\mathcal{Z}(z, 0) = z \sim \rho_0 \in \mathcal{P}^2(M) \\
\rho_T = \mathcal{Z}(\cdot, T){\#} \rho_0 .
\end{cases} 
\end{align}
In the normalizing flow setting, we set the evolution of densities via the geometric pushforward evolution
\begin{align}
\label{eqn:complex_monge_pushfoward}
\log \mathcal{Z}(\cdot,T) \# \rho_0 = \log \rho_0(z) - \log \Big( \frac{\det h(z, \overline{z}, T)}{\det h(z, \overline{z}, 0)} \Big) ,
\end{align}
which is essentially a standard change-of-variables. Assume the hypotheses of the Fundamental Theorem of Calculus are satisfied (the determinant of a Hermitian metric is real-valued, thus we only need the log determinant is absolutely continuous, valid up to the regularity of the Kähler potential). Then
\begin{align}
\log \det h(T) - \log \det h(0) = \int_0^T \frac{\partial}{\partial t}  \log \det h(z, \overline{z}, t)  dt
\end{align}
is the infinitesimal version of the above. Using Jacobi's formula
\begin{align}
\label{eqn:jacobis}
\frac{\partial}{\partial t} \log \det h = \text{Tr}\left( h^{-1} \frac{\partial h}{\partial t} \right) = h^{\overline{j}i} \frac{\partial h_{i\overline{j}}}{\partial t} = h^{\overline{j} i} \partial_i \partial_{\overline{j}} \dot{\Phi} .
\end{align}
Thus \ref{eqn:complex_monge_pushfoward} gives
\begin{align}
\log \mathcal{Z}(\cdot,T) \# \rho_0 = \log \rho_0(z) - \int_0^T \Delta_{\overline{\partial}} \dot{\Phi}  dt .
\end{align}

\subsection{Relations to the instantaneous change of variables}
\label{app:instantaneous_change_of_variables}

\noindent \textbf{Theorem 1 (continued).} \textit{Under the (complex) instantaneous change of variables theorem, the particle vector field obeys
\begin{align}
V = - h^{\overline{j} i} \partial_{\overline{j}} \dot{\Phi} \frac{\partial}{\partial z^i}  ,
\end{align}
up to constants.
}

\vspace{2mm}
\noindent \textit{Proof (slightly informal).} Let $q_t$ denote the probability density of the flow with respect to the evolving metric volume form $\frac{\omega_t^d}{d!}$. The probability measure is therefore $\mu_t = q_t \frac{\omega_t^d}{d!}$. For a normalizing flow governed by a velocity field $V$, mass conservation dictates that the Lie derivative of the measure along the flow must vanish, i.e.
\begin{align}
\frac{\partial}{\partial t} \left( q_t \frac{\omega_t^d}{d!} \right) + \text{div}_{\omega_t} (q_t X) \frac{\omega_t^d}{d!} = 0 .
\end{align}
Denote $X = V + \overline{V}$. Let us expand this via product rule. Using the geometric evolution $\frac{\partial}{\partial t} \left( \frac{\omega_t^d}{d!} \right) = \text{Tr}_{\omega_t}(\dot{\omega}_t) \frac{\omega_t^d}{d!} = (\Delta_{\overline{\partial}} \dot{\Phi}) \frac{\omega_t^d}{d!}$, we obtain the Eulerian continuity equation
\begin{align}
\dot{q}_t + q_t \Delta_{\overline{\partial}} \dot{\Phi} + \text{div}_{\omega_t}(q_t X) = 0 .
\end{align}
Dividing by $q_t$ pointwise (assume sufficient support), we obtain
\begin{align}
\frac{\partial}{\partial t} \log q_t + \Delta_{\overline{\partial}} \dot{\Phi} + \frac{1}{q_t} \text{div}_{\omega_t}(q_t X) = 0 .
\end{align}
Let us choose the $(1,0)$-vector field $V = -(\overline{\partial} \dot{\Phi})^\sharp$, which has a local coordinate representation $V^i = -h^{\overline{j} i} \partial_{\overline{j}} \dot{\Phi}$. The physical divergence can be decomposed into its holomorphic and anti-holomorphic divergences
\begin{align}
\text{div}_{\omega_t}(X) = \nabla_i V^i + \nabla_{\overline{i}} \overline{V}^i \propto \text{Re} \nabla_i V^i .
\end{align}
We will omit the 2 convention for simplicity. By definition
\begin{align}
\text{div}_{\omega_t}(X) \propto \text{Re} \nabla_i V^i = \frac{1}{\det(h)} \frac{\partial}{\partial z^i} \left( \det(h) V^i \right) .
\end{align}
Via product rule
\begin{align}
\text{div}_{\omega_t}(X) \propto \text{Re} \nabla_i V^i = \text{Re} \left( V^i \frac{\partial_i \det(h)}{\det(h)} + \partial_i V^i \right) .
\end{align}
Jacobi's formula gives $\frac{\partial_i \det(h)}{\det(h)} = h^{\overline{m} \ell} \partial_i h_{\ell \overline{m}}$. By the defining property of a Kähler metric, $\partial_i h_{\ell \overline{m}} = \partial_\ell h_{i \overline{m}}$. The derivative of the vector field components is
\begin{align}
\partial_i V^i = \partial_i \big( -h^{\overline{j} i} \partial_{\overline{j}} \dot{\Phi} \big) = -(\partial_i h^{\overline{j} i}) \partial_{\overline{j}} \dot{\Phi} - h^{\overline{j} i} \partial_i \partial_{\overline{j}} \dot{\Phi} .
\end{align}
Using the derivative of the inverse metric $\partial_i h^{\overline{j} i} = -h^{\overline{m} i} (\partial_i h_{\ell \overline{m}}) h^{\overline{j} \ell}$, 
the Christoffel terms cancel, collapsing to the Dolbeault Laplacian
\begin{align}
\text{div}_{\omega_t}(X) = - \text{Re} \ h^{\overline{j} i} \partial_i \partial_{\overline{j}} \dot{\Phi} = -\Delta_{\overline{\partial}} \dot{\Phi} ,
\end{align}
up to constants. The real operator is dropped since $\dot{\Phi}$ is real-valued. Substituting this divergence back into the expanded continuity equation, we notice
\begin{align}
& \frac{\partial}{\partial t} \log q_t + \Delta_{\overline{\partial}} \dot{\Phi} + \langle \nabla \log q_t, X \rangle_{\omega_t} - \Delta_{\overline{\partial}} \dot{\Phi} = \frac{\partial}{\partial t} \log q_t  + \langle \nabla \log q_t, X \rangle_{\omega_t}  = 0 .
\end{align}
Therefore, the material derivative $\frac{d}{dt} \log q_t(Z_t) = 0$. Furthermore, multiplying $\frac{\partial}{\partial t} \log q_t = - \langle \nabla \log q_t, X \rangle_{\omega_t}$ by $q_t \frac{\omega_t^d}{d!}$ recovers the divergence transport $\frac{\partial}{\partial t} \left( q_t \frac{\omega_t^d}{d!} \right) = -\text{div}_{\omega_t} (q_t X) \frac{\omega_t^d}{d!}$.

\section{Relative entropy dissipation}
\label{app:entropy_dissipation}

In this section, we examine the relative entropy (KL divergence) dissipation between the target density and its representation by $\Psi_{k,\theta}$ at a certain time. This is a meaningful way to bridge the statistical qualities of the normalizing flow and the geometric qualities we have been discussing. As we will see in Appendix \ref{app:mabuchi}, the KL divergence has connections to the Mabuchi functional, especially under critical points. We distinguish this from the Mabuchi metric, which is more closely related to an $L^2$ inner product with respect to the form. We note this because we will later discuss the Dirichlet metric.

\vspace{2mm}

\noindent \textbf{Theorem 1 (continued).} \textit{ Denote $p$ the terminal density and $q_t$ the density induced by $\Psi$ at a corresponding step. The first derivative of the KL divergence obeys
\begin{align}
\frac{d}{d t} \text{KL} ( q_t \parallel p) &= - \int_M | \partial \dot{\Phi} |^2 q_t \frac{\omega_t^d}{d!} + \int_M \dot{q}_t \frac{\omega_t^d}{d!} = - \text{Fisher information} + \text{correction} .
\end{align}
Suppose (1) assume $p$ is log-concave $-\sqrt{-1} \partial \overline{\partial} \log p \geq \lambda \omega$ for some $\lambda > 0$; (2) assume $\sqrt{-1} \partial \overline{\partial} \log \frac{q_t}{p} \geq 0$; (3) and the manifold is closed. The second derivative obeys
\begin{align}
\frac{d^2}{dt^2}\text{KL}(q_t\|p) \geq & \left( 2 \lambda   - \frac{C}{8 \nu} \right) \int_M|\partial \dot{\Phi}|^2 q_t \frac{\omega_t^d}{d!}
\\
&  + 2 \text{Re} \int_M \langle \partial \ddot{\Phi}, \partial \dot{\Phi} \rangle_{\omega_t} q_t \frac{\omega_t^d}{d!}
\\ 
& + \int_M \text{Ric}(\partial \dot{\Phi}, \overline{\partial} \dot{\Phi} ) q_t \frac{\omega_t^d}{d!}
\\
& - ( 1 + 2 \nu) \int_M |\nabla^{2,0} \dot{\Phi}|^2 q_t \frac{\omega_t^d}{d!} 
\\
& + \int_M  (\Delta_{\overline{\partial}} \dot{\Phi}) \langle \partial \dot{\Phi}, \partial \log q_t \rangle_{\omega_t} q_t \frac{\omega_t^d}{d!} 
\\
& + \int_M \Big( -(\text{Tr}_{\omega_t}(\mathcal{R}))^2 q_t + h^{\overline{j} \ell} h^{\overline{k} i} \mathcal{R}_{\ell \overline{k}} \mathcal{R}_{i \overline{j}} q_t + h^{\overline{j} i} \partial_t \mathcal{R}_{i \overline{j}} q_t + \text{Tr}_{\omega_t}(\mathcal{R}) \dot{q}_t \Big) \frac{\omega_t^d}{d!}  ,
\end{align}
where $C,\nu$ are real constants independent of dimension, and $\mathcal{R}$ is suitably defined. 
}

\vspace{2mm}

\noindent \textit{Proof (continued in next subsection).} We use notation $\langle \partial f, \partial g \rangle_{\omega_t}$ to denote $h^{\overline{j} i} \partial_i f \partial_{\overline{j}} g$, since the decomposition of forms is orthogonal
\begin{align}
\Lambda^k T^*_{\mathbb{C}}M = \bigoplus_{p+q=k} \Lambda^{p,q} T^*_{\mathbb{C}}M ,
\end{align}
meaning (1,0)-forms are orthogonal to (0,1)-forms, so $\langle \partial f, \overline{\partial} g \rangle_{\omega_t} \equiv 0$. From Jacobi's formula as in \ref{eqn:jacobis}, we have $\frac{\partial}{\partial t}\log(q_t) = \Delta_{\overline{\partial}} \dot{\Phi}$. It is well known the transport of mass satisfies in the continuous case a manifold-variety continuity equation. Exactly, this is
\begin{align}
\dot{q}_t = -\text{div}_{\omega_t}(q_t f) - q_t \text{Tr}_{\omega_t}(\dot{\omega}_t) .
\end{align}
Denote $p$ the terminal density and $q_t$ the density at $t$ according to the normalizing flow. Now, let us consider the KL divergence functional
\begin{align}
\text{KL} ( q_t \parallel p) = \int_M q_t \log ( \frac{q_t}{p} ) \frac{\omega_t^d}{d!}.
\end{align}
Note this KL divergence formulation is different than that used in the loss. The loss functional has two terminal densities, one exact and one the learned pushforward, and is not affected by time, thus the differentiation in time is meaningless with respect to the loss. Differentiating the functional,
\begin{align}
\frac{d}{d t} \text{KL} ( q_t \parallel p) &= \int_M \frac{\partial}{\partial t} \left( q_t \frac{\omega_t^d}{d!} \right) \log q_t - \int_M \frac{\partial}{\partial t} \left( q_t \frac{\omega_t^d}{d!} \right) \log p + \int_M \dot{q}_t \frac{\omega_t^d}{d!} \\
&= -\int_M \text{div}_{\omega_t} (q_t f) \log q_t \frac{\omega_t^d}{d!} + \int_M \text{div}_{\omega_t} (q_t f) \log p \frac{\omega_t^d}{d!} + \int_M \dot{q}_t \frac{\omega_t^d}{d!} ,
\end{align}
and we use the divergence identity $\frac{\partial}{\partial t} \left( q_t \frac{\omega_t^d}{d!} \right) = -\text{div}_{\omega_t} (q_t f) \frac{\omega_t^d}{d!}$. Applying integration by parts with suitable decay, and substituting the gradient flow velocity $f = -(\overline{\partial} \dot{\Phi})^{\sharp}$ which we saw in Appendix \ref{app:instantaneous_change_of_variables} (which leads to a double negative),
\begin{align}
\frac{d}{d t} \text{KL} ( q_t \parallel p) &= -\int_M \langle \partial \dot{\Phi}, \partial \log q_t \rangle_{\omega_t} q_t \frac{\omega_t^d}{d!} + \int_M \langle \partial \dot{\Phi}, \partial \log p \rangle_{\omega_t} q_t \frac{\omega_t^d}{d!} + \int_M \dot{q}_t \frac{\omega_t^d}{d!} 
\\
&= -\int_M \langle \partial \dot{\Phi}, \partial \log \frac{q_t}{p} \rangle_{\omega_t} q_t \frac{\omega_t^d}{d!} + \int_M \dot{q}_t \frac{\omega_t^d}{d!}
\\
&= -\int_M |\partial \dot{\Phi} |^2 q_t \frac{\omega_t^d}{d!} + \int_M \dot{q}_t \frac{\omega_t^d}{d!}.
\end{align}
By definition, $\dot{\Phi} = \log \frac{q_t}{p}$. Therefore, it is also true that
\begin{align}
\dot{F} = \frac{d}{d t} \text{KL}(q_t \parallel p) = -\int_M | \partial  \log \frac{q_t}{p} |^2 q_t  \frac{\omega_t^d}{d!} + \int_M \dot{q}_t \frac{\omega_t^d}{d!} := -I(q_t \parallel p) + \int_M \dot{q}_t \frac{\omega_t^d}{d!} .
\end{align}
Thus, we see
\begin{align}
\frac{d}{dt} \text{KL}(q_t \parallel p ) = -\text{Fisher information} + \text{correction} .
\end{align}
It can be noted from the product rule since by conservation of mass $\frac{\partial}{\partial t} \int_M q_t \frac{\omega_t^d}{d!} = 0$
\begin{align}
\int_M \dot{q}_t \frac{\omega_t^d}{d!} = -\int_M q_t \frac{\partial}{\partial t} \left( \frac{\omega_t^d}{d!} \right) = -\int_M q_t \text{Tr}_{\omega_t}(\dot{\omega}_t) \frac{\omega_t^d}{d!}.
\end{align}

\subsection{Second order normalized relative entropy dissipation}

In this section, we discuss second derivative relative entropy dissipating. We will attempt to show an approximate condition for the second derivative of the KL divergence is nonnegative. This is useful because it establishes a convexity result.

\vspace{2mm}

\noindent Taking the second derivative as in the previous section, and expanding the inner product
\begin{align}
\label{eqn:kl_second_deriv}
\frac{d^2}{dt^2}\text{KL}(q_t \parallel p) = -\frac{\partial}{\partial t}\int_M h^{i\overline{j}}\partial_i\dot{\Phi}\,\partial_{\overline{j}}\dot{\Phi} \cdot q_t \frac{\omega_t^d}{d!} + \frac{\partial}{\partial t} \int_M \dot{q}_t \frac{\omega_t^d}{d!}
\end{align}
We will deal with the second term later. We can note the following evolution equation, which follows under $\partial_t h = \partial_t \text{Ric}$,
\begin{align}
\ddot{\Phi} = - \Delta_{\omega_t} \ddot{\Phi} + | \partial \overline{\partial} \dot{\Phi} |^2 ,
\end{align}
\cite{collins2012twistedkahlerricciflow} the time derivative of the first potential term is
\begin{align} 
-\int_M \langle \partial \ddot{\Phi}, \partial \dot{\Phi}\rangle_{\omega_t} q_t \frac{\omega_t^d}{d!}= \int_M \langle \partial \Delta_{\overline{\partial}}\ddot{\Phi}, \partial \dot{\Phi}\rangle_{\omega_t} q_t \frac{\omega_t^d}{d!} - \int_M \langle \partial | \partial \overline{\partial} \dot{\Phi} |^2, \partial \dot{\Phi} \rangle_{\omega_t} q_t \frac{\omega_t^d}{d!} .
\end{align}
We do not use this identity, since we find the above nontrivial. Thus, we leave this identity purely as a remark. We obtain a second term
\begin{align}
-\int_M \langle \partial \dot{\Phi}, \partial \ddot{\Phi}\rangle_{\omega_t} q_t \frac{\omega_t^d}{d!} .
\end{align}
The time derivative of the inverse metric evolves according to $ \partial_t h^{\overline{j} i} = -h^{\overline{k} i} h^{\overline{j} \ell} \partial_\ell \partial_{\overline{k}} \dot{\Phi} $ (up to a sign convention and an expectation; without loss of generality, take an expectation of the final result), and so the inverse metric term is
\begin{align}
-\int_M h^{\overline{k} i} h^{\overline{j} \ell} \partial_\ell \partial_{\overline{k}} \dot{\Phi} \partial_i\dot{\Phi}\partial_{\overline{j}}\dot{\Phi}\cdot q_t \frac{\omega_t^d}{d!} =  - \int_M \langle \partial \overline{\partial} \dot{\Phi}, \partial \dot{\Phi} \otimes \overline{\partial} \dot{\Phi} \rangle_{\omega_t} q_t \frac{\omega_t^d}{d!} .
\end{align}
Using the identity $\partial_t(\frac{\omega_t^d}{d!}) = \Delta_{\overline{\partial}}\dot{\Phi} \frac{\omega_t^d}{d!} = -\text{Tr}(\mathcal{R}) \frac{\omega_t^d}{d!}$, we arrive at a fourth term
\begin{align}
-\int_M \langle \partial \dot{\Phi}, \partial \log\frac{q_t}{p} \rangle_{\omega_t}(\dot{q}_t - \text{Tr}(\mathcal{R})q_t) \frac{\omega_t^d}{d!}.
\end{align}
Thus, the full expression is
\begin{align}
\label{eqn:second_derivative_KL_full} 
\frac{d^2}{dt^2}\text{KL}(q_t \parallel p) = & \ \int_M \langle \partial \ddot{\Phi}, \partial \dot{\Phi} \rangle_{\omega_t} q_t \frac{\omega_t^d}{d!}
\\
& + \int_M \langle \partial \dot{\Phi}, \partial \ddot{\Phi}\rangle_{\omega_t} q_t  \frac{\omega_t^d}{d!}
\\
& + \int_M \langle \partial \overline{\partial} \dot{\Phi}, \partial \dot{\Phi} \otimes \overline{\partial} \dot{\Phi} \rangle_{\omega_t} q_t \frac{\omega_t^d}{d!}
\\
&-\int_M \langle \partial \dot{\Phi}, \partial \log\frac{q_t}{p}\rangle_{\omega_t}(\dot{q}_t - \text{Tr}(\mathcal{R})q_t) \frac{\omega_t^d}{d!} .
\end{align}
The near repeated inner products can be combined. The first two yield
\begin{align}
\int_M \left( \langle \partial \ddot{\Phi}, \partial \dot{\Phi} \rangle_{\omega_t} + \langle \partial \dot{\Phi}, \partial \ddot{\Phi} \rangle_{\omega_t} \right) q_t \frac{\omega_t^d}{d!} =  2 \text{Re} \int_M \langle \partial \ddot{\Phi}, \partial \dot{\Phi} \rangle_{\omega_t} q_t \frac{\omega_t^d}{d!} .
\end{align}
This follows since
\begin{align}
\langle \partial \ddot{\Phi}, \partial \dot{\Phi} \rangle_{\omega_t} + \langle \partial \dot{\Phi}, \partial \ddot{\Phi} \rangle_{\omega_t} = \langle \partial \ddot{\Phi}, \partial \dot{\Phi} \rangle_{\omega_t} + \overline{\langle \partial \ddot{\Phi}, \partial \dot{\Phi} \rangle}{\omega_t} = 2 \text{Re} \langle \partial \ddot{\Phi}, \partial \dot{\Phi} \rangle_{\omega_t} .
\end{align}
The Laplacian term does not vanish even though the manifold is closed due to $q_t$. For the fourth term, notice $\dot{q}_t$ satisfies a continuity equation, and so
\begin{align}
-\int_M \langle \partial \dot{\Phi}, \partial \log \frac{q_t}{p} \rangle_{\omega_t} (\dot{q}_t - \text{Tr}(\mathcal{R}) q_t)  \frac{\omega_t^d}{d!}= -\int_M \langle \partial \dot{\Phi}, \partial \log \frac{q_t}{p} \rangle_{\omega_t} (-\text{div}_{\omega_t}(q_t (-\overline{\partial} \dot{\Phi})^{\sharp}))  \frac{\omega_t^d}{d!}.
\end{align}
We remark the above uses a triple negative. We must include the $\mathcal{R}$ term into the divergence due to $\frac{\partial}{\partial t} \left( q_t \frac{\omega_t^d}{d!} \right) = -\text{div}_{\omega_t} (q_t f) \frac{\omega_t^d}{d!}$. Shifting the divergence,
\begin{align}
-\int_M \langle \partial\dot{\Phi}, \partial \dot{\Phi} \rangle_{\omega_t} \text{div}_{\omega_t}(q_t (\overline{\partial}\dot{\Phi})^{\sharp})  \frac{\omega_t^d}{d!}=  \int_M \langle \partial \langle \partial \dot{\Phi}, \partial \dot{\Phi} \rangle, \partial \dot{\Phi} \rangle_{\omega_t} q_t  \frac{\omega_t^d}{d!}.
\end{align}
Expanding the outer gradient using product rule for connections and using Hessian notation,
\begin{align}
& \int_M \left[ \langle \nabla_{\partial \dot{\Phi}} \partial \dot{\Phi}, \partial \dot{\Phi} \rangle_{\omega_t} + \langle \partial \dot{\Phi}, \nabla_{\partial \dot{\Phi}} \partial \dot{\Phi} \rangle_{\omega_t} \right] q_t \frac{\omega_t^d}{d!} 
\\
& = \int_M \left[ \nabla^{2,0} \dot{\Phi} (\overline{\partial} \dot{\Phi}, \overline{\partial} \log \frac{q_t}{p}) + \partial \overline{\partial} \dot{\Phi} (\partial \dot{\Phi}, \overline{\partial} \dot{\Phi}) \right] q_t  \frac{\omega_t^d}{d!} .
\end{align}
Since $\nabla^{2,0} \dot{\Phi}$ is a (2,0)-Hessian as a bilinear form, it takes two (0,1)-form input. Putting everything together, we conclude
\begin{align}
\label{eqn:1}
\frac{d^2}{dt^2} \text{KL}(q_t \parallel p) = &  \ 2 \text{Re} \int_M \langle \partial \ddot{\Phi}, \partial \dot{\Phi} \rangle_{\omega_t} q_t \frac{\omega_t^d}{d!}
\\
& + \int_M  \nabla^{2,0} \dot{\Phi} (\overline{\partial} \dot{\Phi}, \overline{\partial} \log \frac{q_t}{p}) q_t \frac{\omega_t^d}{d!}
\\
& + 2 \int_M \langle \partial \overline{\partial} \dot{\Phi}, \partial \dot{\Phi} \otimes \overline{\partial} \dot{\Phi} \rangle_{\omega_t} q_t \frac{\omega_t^d}{d!} .
\end{align}
Let us bound the Hessian term, which is nontrivial. Let us isolate the $\log q_t$, hence $q_t \overline{\partial} \log q_t = \overline{\partial} q_t$. Notice
\begin{align}
\int_M \nabla^{2,0} \dot{\Phi} (\overline{\partial} \dot{\Phi}, \overline{\partial} \log q_t) q_t \frac{\omega_t^d}{d!} & = \int_M h^{\overline{k} i} h^{\overline{l} j} \nabla_i \nabla_j \dot{\Phi} \nabla_{\overline{k}} \dot{\Phi} \nabla_{\overline{l}} q_t \frac{\omega_t^d}{d!} \\
& = - \int_M \nabla_{\overline{l}} \left( h^{\overline{k} i} h^{\overline{l} j} \nabla_i \nabla_j \dot{\Phi} \nabla_{\overline{k}} \dot{\Phi} \right) q_t \frac{\omega_t^d}{d!} \\
& = - \int_M h^{\overline{k} i} h^{\overline{l} j} \left( \nabla_{\overline{l}} \nabla_i \nabla_j \dot{\Phi} \nabla_{\overline{k}} \dot{\Phi} + \nabla_i \nabla_j \dot{\Phi} \nabla_{\overline{l}} \nabla_{\overline{k}} \dot{\Phi} \right) q_t \frac{\omega_t^d}{d!}.
\end{align} 
We have noted the metric is compatible in the last line. In particular, we get
\begin{align}
\int_M \nabla^{2,0} \dot{\Phi} (\overline{\partial} \dot{\Phi}, \overline{\partial} \log q_t) q_t \frac{\omega_t^d}{d!} = - \int_M (\nabla^j \nabla_i \nabla_j \dot{\Phi}) \nabla^i \dot{\Phi} q_t \frac{\omega_t^d}{d!} - \int_M |\nabla^{2,0} \dot{\Phi}|^2 q_t \frac{\omega_t^d}{d!} .
\end{align}
Returning to the original term,
\begin{align}
\int_M \nabla^{2,0} \dot{\Phi} & (\overline{\partial} \dot{\Phi}, \overline{\partial} \log \frac{q_t}{p}) q_t \frac{\omega_t^d}{d!} = \int_M \nabla^{2,0} \dot{\Phi} (\overline{\partial} \dot{\Phi}, \overline{\partial} \log q_t) q_t \frac{\omega_t^d}{d!} - \int_M \nabla^{2,0} \dot{\Phi} (\overline{\partial} \dot{\Phi}, \overline{\partial} \log p) q_t \frac{\omega_t^d}{d!} 
\\
& = - \int_M |\nabla^{2,0} \dot{\Phi}|^2 q_t \frac{\omega_t^d}{d!} - \int_M (\nabla^j \nabla_i \nabla_j \dot{\Phi}) \nabla^i \dot{\Phi} q_t \frac{\omega_t^d}{d!} - \int_M \nabla^{2,0} \dot{\Phi} (\overline{\partial} \dot{\Phi}, \overline{\partial} \log p) q_t \frac{\omega_t^d}{d!} 
\\
& \geq - \int_M |\nabla^{2,0} \dot{\Phi}|^2 q_t \frac{\omega_t^d}{d!} - \int_M (\nabla^j \nabla_i \nabla_j \dot{\Phi}) \nabla^i \dot{\Phi} q_t \frac{\omega_t^d}{d!} 
\\
& \ \ \ \ \ \ \ \ \ \ \ \ \ \ - 2\nu \int_M |\nabla^{2,0} \dot{\Phi}|^2 q_t \frac{\omega_t^d}{d!} - \frac{1}{8\nu} \int_M |\overline{\partial} \dot{\Phi}|^2 |\overline{\partial} \log p|^2 q_t \frac{\omega_t^d}{d!} 
\\
& \geq - (1+2\nu) \int_M |\nabla^{2,0} \dot{\Phi}|^2 q_t \frac{\omega_t^d}{d!} - \int_M (\nabla^j \nabla_i \nabla_j \dot{\Phi}) \nabla^i \dot{\Phi} q_t \frac{\omega_t^d}{d!} - \frac{C}{8\nu} \int_M |\overline{\partial} \dot{\Phi}|^2 q_t \frac{\omega_t^d}{d!}. 
\end{align}
The first inequality is the Young's (Peter-Paul) inequality. Let us examine the middle term in the last line of the above. Observe the identity $\nabla^j \nabla_i \nabla_j \dot{\Phi} = \nabla_i \Delta_{\overline{\partial}} \dot{\Phi} - \text{Ric}_i^\ell \nabla_\ell \dot{\Phi}$. Thus, and by integrating by parts,
\begin{align}
-\int_M (\nabla^j \nabla_i \nabla_j \dot{\Phi}) & \nabla^i \dot{\Phi} q_t \frac{\omega_t^d}{d!} = -\int_M \langle \partial \Delta_{\overline{\partial}} \dot{\Phi}, \partial \dot{\Phi} \rangle_{\omega_t} q_t \frac{\omega_t^d}{d!} + \int_M \text{Ric}(\partial \dot{\Phi}, \overline{\partial} \dot{\Phi}) q_t \frac{\omega_t^d}{d!}
\\
& = \int_M (\Delta_{\overline{\partial}} \dot{\Phi})^2 q_t \frac{\omega_t^d}{d!} + \int_M (\Delta_{\overline{\partial}} \dot{\Phi}) \langle \partial \dot{\Phi}, \partial \log q_t \rangle_{\omega_t} q_t \frac{\omega_t^d}{d!} + \int_M \text{Ric}(\partial \dot{\Phi}, \overline{\partial} \dot{\Phi}) q_t \frac{\omega_t^d}{d!}
\\
& \geq \int_M  (\Delta_{\overline{\partial}} \dot{\Phi}) \langle \partial \dot{\Phi}, \partial \log q_t \rangle_{\omega_t} q_t \frac{\omega_t^d}{d!} + \int_M \text{Ric}(\partial \dot{\Phi}, \overline{\partial} \dot{\Phi}) q_t \frac{\omega_t^d}{d!}.
\end{align}
Turning to the last term in \ref{eqn:1}, under the log-concavity assumption,
\begin{align}
\int_M \langle \partial \overline{\partial} \dot{\Phi}, \partial \dot{\Phi} \otimes \overline{\partial} \dot{\Phi} \rangle_{\omega_t} q_t  \frac{\omega_t^d}{d!}\geq \lambda \int_M |\partial \dot{\Phi}|^2 q_t  \frac{\omega_t^d}{d!}.
\end{align}
Now, let us return to the second term as in \ref{eqn:kl_second_deriv}. Notice
\begin{align}
\int_M \dot{q}_t \frac{\omega_t^d}{d!} = -\int_M q_t \frac{\partial}{\partial t}\left(\frac{\omega_t^d}{d!}\right) = \int_M \text{Tr}_{\omega_t}(\mathcal{R}) q_t \frac{\omega_t^d}{d!} .
\end{align}
\textit{Claim.} We have
\begin{align}
\frac{d}{dt} \int_M \text{Tr}_{\omega_t}(\mathcal{R}) q_t \frac{\omega_t^d}{d!} = \int_M \Big( -(\text{Tr}_{\omega_t}(\mathcal{R}))^2 q_t + h^{\overline{j} \ell} h^{\overline{k} i} \mathcal{R}_{\ell \overline{k}} \mathcal{R}_{i \overline{j}} q_t + h^{\overline{j} i} \partial_t \mathcal{R}_{i \overline{j}} q_t + \text{Tr}_{\omega_t}(\mathcal{R}) \dot{q}_t \Big) \frac{\omega_t^d}{d!} ,
\end{align}
where $\mathcal{R}$ is a custom tensor to match what we established in \ref{app:statistical_densities}.

\vspace{2mm}

\noindent \textit{Proof of claim.} Observe
\begin{align}\frac{d}{d t} \int_M \text{Tr}_{\omega_t}(\mathcal{R}) q_t \frac{\omega_t^d}{d!} = \int_M \frac{\partial (\text{Tr}_{\omega_t}(\mathcal{R}))}{\partial t} q_t \frac{\omega_t^d}{d!} + \int_M \text{Tr}_{\omega_t}(\mathcal{R}) \dot{q}_t \frac{\omega_t^d}{d!} + \int_M \text{Tr}_{\omega_t}(\mathcal{R}) q_t \frac{\partial}{\partial t}\left( \frac{\omega_t^d}{d!} \right) .
\end{align}
Substituting the volume form derivative $\frac{\partial}{\partial t}\left( \frac{\omega_t^d}{d!} \right) = -\text{Tr}_{\omega_t}(\mathcal{R}) \frac{\omega_t^d}{d!}$,
\begin{align}\frac{d}{d t} \int_M \text{Tr}_{\omega_t}(\mathcal{R}) q_t \frac{\omega_t^d}{d!} = \int_M \frac{\partial (\text{Tr}_{\omega_t}(\mathcal{R}))}{\partial t} q_t \frac{\omega_t^d}{d!} + \int_M \text{Tr}_{\omega_t}(\mathcal{R}) \dot{q}_t \frac{\omega_t^d}{d!} - \int_M (\text{Tr}_{\omega_t}(\mathcal{R}))^2 q_t \frac{\omega_t^d}{d!} .
\end{align}
Under the statistical evolution, the trace of the effective curvature evolves as
\begin{align}
\frac{\partial (\text{Tr}_{\omega_t}(\mathcal{R}))}{\partial t} = h^{\overline{j} \ell} h^{\overline{k} i} \mathcal{R}_{\ell \overline{k}} \mathcal{R}_{i \overline{j}} + h^{\overline{j} i} \partial_t \mathcal{R}_{i \overline{j}} .
\end{align}
Substituting this directly into the integral proves the claim.

\noindent $ \square $

\vspace{2mm}

\noindent Putting everything together, and noting $|\partial \dot{\Phi}|^2 = |\overline{\partial} \dot{\Phi}|^2 = h^{\overline{j} i}\partial_i \dot{\Phi} \partial_{\overline{j}}\dot{\Phi}$, we have the result. Now, we can also note
\begin{align}
\ddot{\Phi} = \frac{\dot{q}_t}{q_t}  .
\end{align}
This section concludes Theorem 1.

\noindent $ \square $

\section{Relations to other topics in geometry}

\subsection{Kähler-Einstein conditions}
\label{app:kahler_einstein_conditions}

The main results and objectives have structural similarities to the traditional Kähler-Einstein condition
\begin{align}
\lambda h_{i \overline{j}} = \text{Ric}_{i \overline{j}} ,
\end{align}
and so a lot of our work is reminiscent of this condition here. Recall in section \ref{sec:training}, we were interested in the pointwise differentiation of the loss  $\partial_i \partial_{\overline{j}} | \Psi_\theta^{-1} |^2 -  \text{Ric}_{i \overline{j}}$. Recall $\partial_i \partial_{\overline{j}} |\Psi_{\theta}^{-1}|^2$ corresponds to the unit Gaussian scenario. Generalized and with Ricci curvature, this corresponds to $-\partial_i \partial_{\overline{j}} \log p_{\alpha}(\Psi_{\theta}^{-1}(w)) -  \text{Ric}_{i \overline{j}}$, where the first term is now a Fisher metric. We would have at unit Gaussian
\begin{align}
\partial_i \partial_{\overline{j}} | \Psi_\theta^{-1} |^2 = \delta_{i \overline{j}} = h_{i \overline{j}} =  \text{Ric}_{i \overline{j}} .
\end{align}
The difference has become a singular Ricci curvature term because the Ricci curvature of the isotropic unit Gaussian is zero, i.e. $\log \det(h) = 0$ when $h = \delta$, which is a Ricci-flat scenario. In the general case,
\begin{align}
-\EX \partial_i \partial_{\overline{j}} \log p_{k,\alpha}(\Psi_{\theta}^{-1}(w)) - \EX   \text{Ric}_{k,i \overline{j}} +\EX \partial_i \partial_{\overline{j}} \log p_{k-1,\alpha}(\Psi_{\theta}^{-1}(w)) + \EX   \text{Ric}_{k-1,i \overline{j}}  = 0 ,
\end{align}
which have similarity to Kähler-Einstein conditions. In particular, if the two terms split and are individually zero then it is Kähler-Einstein. From our main result in \ref{sec:main_contribution}, we have
\begin{align}
\partial_t h + \partial_t \EX \text{Ric} ,
\end{align}
which is almost Kähler-Einstein. The Kähler-Einstein condition is a special case of the Kähler-Ricci soliton equation $\text{Ric}_{i \overline{j}} = \lambda h_{i \overline{j}}$ when $X \in \Gamma(T^{(1,0)}M)$ is a Killing field in the holomorphic tangent bundle, and so the Lie derivative term $\mathcal{L}_X h_{i \overline{j}}$ vanishes, i.e.
\begin{align}
\text{Ric}_{i \overline{j}} - \lambda h_{i \overline{j}} = \mathcal{L}_X h_{i \overline{j}} = \mathcal{L}_X \text{Ric}_{i \overline{j}} = X^k \partial_k \text{Ric}_{i \overline{j}} + \overline{X}^k \partial_{\overline{k}} \text{Ric}_{i \overline{j}} + \text{Ric}_{k \overline{j}} \partial_i X^k + \text{Ric}_{i \overline{k}} \partial_{\overline{j}} \overline{X}^k = 0 .
\end{align}
The second to last equality is an expansion of the Lie derivative. The constant, which is $\lambda$ in this scenario, is the Einstein or cosmological constant. We remark normalizing flow loss minimization is typically not at a zero value, although recall the failure of dominated convergence in \ref{sec:training}. Also, this scenario would correspond to the first Chern class as zero in the de Rham cohomology sense, so this complements the next subsection \ref{app:calabi-yau}.

\subsection{Relations to first Chern class and the Calabi-Yau manifold}
\label{app:calabi-yau}

Let $h(t) \in \Gamma(M, T^{*(1,0)} M \otimes T^{*(0,1)}M)$ be a Kähler metric on $M$. Recall the Dolbeault operators $(\partial, \overline{\partial})$ can act on $(0,0)$-forms, or scalar-valued functions. Since the Ricci form represents $2\pi c_1(M)$ in the de Rham cohomology, one has
\begin{align}
[-\sqrt{-1} \partial \overline{\partial} \log \det  h ] = 2 \pi c_1(M) ,
\end{align}
where $[\cdot]$ is the cohomology class and $c_1$ is the first Chern class. In the case that $c_1(M)=0$, it follows that $(M,h)$ is consistent with a Calabi-Yau manifold \cite{chu2024Kahlermanifoldsnonnegativemixed}. In this case, there is some correspondence with the Calabi-Yau volume form identity
\begin{align}
\frac{\omega_t^d}{d!} = \left( \frac{\sqrt{-1}}{2} \right)^d \Omega \wedge \overline{\Omega} \cdot (\text{constant}) ,
\end{align}
where $\Omega$ is a nowhere-vanishing holomorphic $d$-form, which relates the volume form we have been working with to the Calabi-Yau manifold.

\section{Relations to vector fields}

\subsection{Inducing $\Psi$ from $f$ and $\Phi$}
\label{app:sharp_section}

Let us examine the anti-holomorphic differential form in the cotangent bundle
\begin{align}
\alpha_t = \sum_i \frac{\partial \dot{\Phi}_t}{\partial \overline{z}^i} d\overline{z}^i = \overline{\partial} \dot{\Phi}_t \in \Gamma(T^{*(0,1)}M) .
\end{align}
We noted in section \ref{app:instantaneous_change_of_variables} that $f$ has a Kähler formulation, which has an equivalent musical isomorphism formulation
\begin{align}
\label{eqn:f_phi_sharp}
f_t = -(\overline{\partial} \dot{\Phi}_t)^\sharp = -h^{\overline{j} i} \partial_{\overline{j}} \dot{\Phi}_t \frac{\partial}{\partial z^i} ,
\end{align}
where $\sharp : \Gamma(T^{*(0,1)}M) \rightarrow \Gamma(T^{(1,0)}M)$ is the Hermitian sharp musical isomorphism operator. Therefore, integrating the ODE and substituting in \ref{eqn:f_phi_sharp}, we arrive at the system
\begin{align}
\begin{cases}
& \frac{d}{d t} \Psi_{t,\theta}(z_0) = f(\Psi_{t,\theta}(z_0), t) 
\\
& \Psi_{t,\theta}(z_0) = z_0 - \int_0^t (\overline{\partial} \dot{\Phi}_s)^{\sharp} (\Psi_{s,\theta}(z_0)) ds  ,
\end{cases}
\end{align}
and so $\dot{\Psi}_{\theta}^i(z) = -h^{\overline{j} i} \partial_{\overline{j}} \dot{\Phi}$ by differentiating and using the sharp operator. We will use this in Appendix \ref{app:perelman_functionals}. In particular, we have found a way to relate
\begin{align}
\alpha_t \in \Gamma(T^{*(0,1)}M) \xrightarrow{\sharp} f_t \in \Gamma(T^{(1,0)}M) ,
\end{align}
or with the flat map $\flat$, $f_t^{\flat} = - \overline{\partial} \dot{\Phi}_t$.

\subsection{Relations to kinetic energies}
\label{app:kinetic_energies}

A trajectory of data following the normalizing flow in complex space $\mathbb{C}^d$ can be integrated against to measure the amount of work done by the data particle. In particular, we examine
\begin{align}
\int_{\gamma} \langle v, (\overline{\partial} \dot{\Phi})^\sharp \rangle_{\omega_t} dt,
\end{align}
where $\gamma$ is a curve in complex space and $v$ is the velocity vector. Using Kähler potential $\Phi$, the above integral uses the velocity and rewrites it as
\begin{align}
\int_0^T h^{\overline{j} i}  \partial_{i} \dot{\Phi}  \partial_{\overline{j}} \dot{\Phi}  dt ,
\end{align}
which is precisely a kinetic energy functional. We show this more rigorously. From Appendix \ref{app:sharp_section}, we saw $\dot{\Psi}_{\theta} = -(\overline{\partial} \dot{\Phi})^\sharp$. Thus, since $\Psi_{\theta}$ is non-holomorphic, so is $\dot{\Psi}_{\theta}$, and it must be true the (0,1) part does not vanish under the $\overline{\partial}$ operator, with slight notation abuse,
\begin{align}
\overline{\partial}(\overline{\partial} \dot{\Phi})^\sharp = \overline{\partial} ( z_0 + \int_0^t (\overline{\partial} \dot{\Phi})^{\sharp} dt) \neq 0.
\end{align}
In particular, the chain is 
\begin{align}
\dot{\Phi}_t \in  C^{\infty}(M) \stackrel{\overline{\partial}} {\rightarrow} \overline{\partial} \dot{\Phi}_t \in \Omega^{0,1}(M) \stackrel{\sharp}{\rightarrow} (\overline{\partial} \dot{\Phi}_t)^{\sharp} \in \Gamma(M,T^{(1,0)}M) \stackrel{\overline{\partial}}{\not\rightarrow} 0 .
\end{align}
Using the definition of kinetic energy,
\begin{align}
\text{kinetic energy} := \int_0^T |\dot{\Psi}_{\theta}|^2_{\omega_t}  dt = \int_0^T |(\overline{\partial} \dot{\Phi})^{\sharp}|^2_{\omega_t}  dt = \int_0^T h^{\overline{j}i} \partial_i \dot{\Phi}  \partial_{\overline{j}} \dot{\Phi}  dt \geq 0.
\end{align}
Using $v = \dot{\Psi}_{\theta}$, we get the action
\begin{align}
\text{action} = \int_{\gamma} \langle v, (\overline{\partial} \dot{\Phi})^\sharp \rangle_{\omega_t} dt = \int_{0}^{T} \langle -(\overline{\partial} \dot{\Phi})^\sharp, (\overline{\partial} \dot{\Phi})^\sharp \rangle_{\omega_t} dt = - \int_{0}^{T} |v|^2_{\omega_t} dt .
\end{align}
Note that this is closely related to the geometric optimal transport problem as in Appendix \ref{app:continuity_equations}. The above kinetic energy functional has close connections to the Mabuchi metric and scalar curvature. If the particle's curve is closed, then applying complex Stokes' theorem and under the (unnormalized) Kähler-Ricci flow
\begin{align}
\label{eqn:stokes}
\oint_{\gamma} \partial \dot{\Phi} = \iint_{D} \overline{\partial} \partial \dot{\Phi} = - \iint_{D} \partial \overline{\partial} \dot{\Phi} = \sqrt{-1} \iint_D dd^c \dot{\Phi} = -\sqrt{-1} \iint_D \text{Ric}(\omega)  .
\end{align}
Here, $\partial \dot{\Phi} = \sum_{i=1}^d \frac{\partial \dot{\Phi}}{\partial z^i} dz^i$ is a (1,0)-form. We crucially remark Cauchy's residue does not necessarily apply because $\partial \dot{\Phi}$ is not a meromorphic $(1,0)$-form, but under more relaxed conditions we notice
\begin{align}
\oint_{\gamma} \partial \dot{\Phi} = 2\pi \sqrt{-1} \sum_{z = p_k} \text{Res}(\partial \dot{\Phi}) .
\end{align}
Here, $d = \partial + \overline{\partial}$ and $d^c = -\frac{\sqrt{-1}}{2}(\partial - \overline{\partial})$ are the convention so that $dd^c = \sqrt{-1} \partial \overline{\partial}$, thus the third equality is purely notational. The last equality in \ref{eqn:stokes} follows by the definition of the form evolution under the Kähler-Ricci flow.

\section{Functional critical point equivalences}
\label{app:mabuchi}

In this section, we provide a result that the first variation using the KL divergence matches that of the Mabuchi functional. It is well known that the critical points of the Mabuchi functional correspond to Kähler-Einstein geometries. 

\vspace{2mm}

\noindent In this section, we will assume \cite{he2025twistedcalabifunctionaltwisted}
\begin{align}
\text{Tr}_{\omega_t}(\rho_K) = \overline{R}
\end{align}
pointwise, which is a twisted Kähler-Einstein condition. We will define $\rho_K$ later. In general, the above does not hold. Therefore, the two functionals have only the same critical points.

\vspace{2mm}

\noindent \textbf{Theorem 2 (continued).} \textit{ If $\mathcal{M}$ is the Mabuchi functional, then it has the same critical points as the KL divergence at a Kähler scalar curvature condition.
}

\vspace{2mm}

\noindent \textit{Proof.} Let us denote $\omega_0^d$ a reference volume form, $q$ the target density, and $f$ the pushforward of the initial density. Let us consider the KL divergence as a functional
\begin{align}
\text{KL}(f \parallel q) = \int_M f \log ( \frac{f}{q} ) \frac{\omega_0^d}{d!} .
\end{align}
Notice this definition is inverted from that we used in training. We can note $f = \frac{\omega_t^d}{\omega_0^d}$. The first variation of the KL divergence can be derived as
\begin{align}
\frac{d}{d \epsilon} |_{\epsilon = 0} \text{KL}( f + \epsilon \delta f \parallel q ) = \frac{d}{d \epsilon} |_{\epsilon = 0} \int_M ( f + \epsilon \delta f) \log (\frac{ f + \epsilon \delta f}{q} ) \frac{\omega_0^d}{d!} ,
\end{align}
and is given by
\begin{align}
\delta \mathcal{F} = \int_M \delta f \left[  \log ( \frac{f}{q} ) + 1 \right] \frac{\omega_0^d}{d!} = \int_M \log \left( \frac{f}{q} \right) \frac{\delta \omega_t^d}{d!} .
\end{align}
The constant 1 is eliminated upon integration since $\delta f$ is signed to preserve the probability mass condition.  Recall the identity for the first variation of the volume form \cite{calabi2001spacekahlermetricsii} \cite{collins2018inversemongeampereflowapplications}
\begin{align}
\delta \omega_t^d = (\Delta_{\overline{\partial}} \delta \Phi) \omega_t^d .
\end{align}
Substituting into the KL variation,
\begin{align}
\delta \mathcal{F} = \int_M \log( \frac{f}{q}) \Delta_{\overline{\partial}} \delta \Phi \frac{\omega_t^d}{d!} .
\end{align}
Applying Green's second identity with a decay condition, i.e. the Laplacian is self-adjoint,
\begin{align}
\delta \mathcal{F} & = \int_M \delta \Phi \Delta_{\overline{\partial}} \log ( \frac{ f}{ q } ) \frac{\omega_t^d}{d!} 
\\
 & = \int_M \delta \Phi  [ \Delta_{\overline{\partial}} \log (f) - \Delta_{\overline{\partial}} \log(q) ] \frac{\omega_t^d}{d!}
 \\
& \stackrel{(1)}{=}  \int_M \delta \Phi [ \text{Tr}_{\omega_t} \text{Ric}(\omega_0^d) - R + (\text{Tr}_{\omega_t}\rho_K - \text{Tr}_{\omega_t} \text{Ric}(\omega_0^d)  )] \frac{\omega_t^d}{d!} 
\\
& = - \int_M \delta \Phi [ R - \text{Tr}_{\omega_t} ( \rho_K) ] \frac{\omega_t^d}{d!}.
\end{align}
(1) follows since
\begin{align}
f = \frac{\omega_t^d}{\omega_0^d} .
\end{align}
Therefore, the difference in Ricci forms satisfies
\begin{align}
\rho_t - \rho_0 = -\sqrt{-1} \partial \overline{\partial} \log \left( \frac{\omega_t^d}{\omega_0^d} \right) = -\sqrt{-1} \partial \overline{\partial} \log f ,
\end{align}
and so
\begin{align}
\Delta_{\overline{\partial}} \log f = \text{Tr}_{\omega_t} (\sqrt{-1} \partial \overline{\partial} \log f) .
\end{align}
Since scalar curvature is the trace of $\rho_t$, i.e. $\text{Tr}_{\omega_t} \rho_t = R$, we get
\begin{align}
\Delta_{\overline{\partial}} \log f = \text{Tr}_{\omega_t} \text{Ric}(\omega_0) - R .
\end{align}
We have denoted the twisted Ricci form $\rho_K = - \sqrt{-1} \partial \overline{\partial} \log (q) + \text{Ric}(\omega_0)$ \cite{george2025complexmongeampereequationpositive}. Taking the trace, $-\Delta_{\overline{\partial}} \log q = \text{Tr}_{\omega_t} \rho_K - \text{Tr}_{\omega_t} \text{Ric}(\omega_0)$. This is equivalent to the Mabuchi metric up to average scalar curvature, i.e. 
\begin{align}
\delta \mathcal{M} = - \int_M \delta \Phi ( R - \overline{R}) \frac{\omega_t^d}{d!} .
\end{align}
Thus the loss possesses the identity
\begin{align}
\delta \mathcal{F} =   \delta \mathcal{M} .
\end{align}
By the initial pointwise assumption, the result follows. As a side remark, it can be noted
\begin{align}
\text{Tr}_{\omega} (\alpha) \frac{\omega^d}{d!} = \alpha \wedge  \frac{\omega^{d-1}}{(d-1)!} .
\end{align}
Thus,
\begin{align}
\int_M \text{Tr}_{\omega_t} (\rho_K) \frac{\omega_t^d}{d!} = \int_M \rho_K \wedge \frac{\omega_t^{d-1}}{(d-1)!} .
\end{align}
We can see
\begin{align}
\int_{M} \text{Tr}_{\omega_{t}}(\rho_{K}) \frac{\omega_t^d}{d!} &= \int_{M} \frac{1}{(d-1)!} \rho_{K} \wedge \omega_{t}^{d-1} \\
&= \int_{M} \frac{1}{(d-1)!} (\text{Ric}(\omega_{t}) + \sqrt{-1} \partial\overline{\partial}\psi) \wedge \omega_{t}^{d-1} \\
&\stackrel{\text{Stokes' theorem}}{=} \int_{M} \frac{1}{(d-1)!} \text{Ric}(\omega_{t}) \wedge \omega_{t}^{d-1} + \int_{M} \frac{1}{(d-1)!} d(\sqrt{-1} \ \overline{\partial}\psi \wedge \omega_{t}^{d-1}) \\
&= \int_{M} R \frac{\omega_t^d}{d!} + 0 \\
&= V \cdot \overline{R} .
\end{align}
We have used $d \omega_t^d = 0$.

\section{Kähler-Ricci relations to Perelman-type functionals under Dirichlet metrics}
\label{app:perelman_functionals}

In this section, we study gradient flows for the unnormalized Kähler-Ricci flow. As we saw in section \ref{app:normalized_kahler_ricci_flow}, the complex normalizing flow associated to the normalized equation is less standard from a machine learning perspective and is of overall less interest to us. We will not study the scalar curvature variety of this functional but rather a form variety. The reason for this is we can establish a gradient flow on the potential otherwise easily, relating the log of the form quotient to the flow map. Noting a relationship with the Dirichlet gradient and the Kähler potential, we may establish a relation between $\dot{\Phi}$ and $\nabla \Psi$ up to a function on the manifold $f$ \cite{perelman2002entropyformularicciflow}. Ricci flow is not a gradient flow in the space of all metrics. Perelman showed Ricci flow is a gradient flow with respect to $f$ of a fixed measure, thus $f$ acts as a weighted measure so that the metric evolves in accordance with $f$.

\vspace{2mm}

\noindent In the normalized case, it is known that the Kähler-Ricci flow is a gradient flow of the Ding functional \cite{collins2018inversemongeampereflowapplications}
\begin{align}
\mathcal{D}(\Phi) = -\mathcal{E}(\Phi) - \log \Bigg( \int_M e^{f - \Phi} \omega^d \Bigg) ,
\end{align}
where $\mathcal{E}(\Phi)$ is the Aubin-Yau functional, whose variation satisfies
\begin{align}
\delta \mathcal{E} = \int_M \delta \Phi \frac{\omega_\Phi^d}{d!} .
\end{align}
Note the critical points correspond to Kähler-Einstein metrics, thus we also refer to section \ref{app:kahler_einstein_conditions}. In this section, we will use notation $\omega_{\Phi}^d$ instead of $\omega_t^d$ to clarify when a functional perturbation is taken, i.e. we will use $\omega_{\Phi + \epsilon \delta \Phi}^d$.

\vspace{2mm}

\noindent \textbf{Theorem 2 (continued).} \textit{ Under Kähler-Ricci flow, suitable $f$ and the Dirichlet metric, $\Phi$ is governed by a gradient flow satisfying
\begin{align}
\dot{\Phi} = -\doublenabla \mathcal{F}_K =  \log \det h - f  .
\end{align}
}

\vspace{2mm}

 \noindent \textit{Proof.} Perelman's functional is classically
\begin{align}
\label{eqn:perelman_func}
\mathcal{F}(h, f) = \int_M (R + |\nabla f|^2) e^{-f} \frac{\omega_t^d}{d!}.
\end{align}
We examine a Perelman-type functional but using the Kähler potential as opposed to a scalar curvature formulation. In the unnormalized case, we examine the Perelman-type functional \cite{klemyatin2025convergenceinversemongeampereflow} \cite{shen2022cherncalabiflowhermitianmanifolds}
\begin{align}
\mathcal{F}_{K}(\Phi) = \int_M \left( \log \frac{\omega_\Phi^d}{\omega^d} - f \right) \frac{\omega_\Phi^d}{d!} ,
\end{align}
which is mostly a form-log version of \ref{eqn:perelman_func}. Here $\omega$ is the base form, i.e. that of the base distribution or complex unit Gaussian, and $\omega_{\Phi}$ is the transformed form, i.e. that pertaining to $\Psi_{k,\theta}$ for some $k$ in the normalizing flow. $f$ is a Ricci potential of the base, and has connections to Ricci flow under diffeomorphisms \cite{perelman2002entropyformularicciflow}. Under the Dirichlet metric, we get
\begin{align}
\doublenabla \mathcal{F}_K = \log \frac{\omega_\Phi^d}{\omega^d} - f  .
\end{align}
It is the unique tangent vector satisfying the Riesz representation property $\delta \mathcal{F}_K(v) = \langle \doublenabla \mathcal{F}_K, v \rangle_\Phi$, and $\delta \mathcal{F}_K$ is the first variation (in other words, $\doublenabla \mathcal{F}_K$ is the functional derivative). Thus the gradient flow recovers \cite{Cao1985} \cite{phong2007multiplieridealsheaveskahlerricci}
\begin{align}
\dot{\Phi} = \log \frac{\omega_\Phi^d}{\omega^d} - f .
\end{align}
This has correspondence to an unnormalized Kähler-Ricci flow
\begin{align}
\frac{\partial \omega_\Phi}{\partial t} = -\text{Ric}(\omega_\Phi) .
\end{align}
Let us verify this, specializing to our normalizing flows. We will take the first variation. First, we perturb
\begin{align}
A(\Phi + \epsilon\,\delta\Phi) = \int_M \left(\log\frac{\omega_{\Phi+\epsilon\,\delta\Phi}^d}{\omega^d} - f\right)\frac{\omega_{\Phi+\epsilon\,\delta\Phi}^d}{d!} .
\end{align}
Given the form variation $\omega_{\Phi+\epsilon \delta\Phi} = \omega_\Phi + \epsilon \sqrt{-1}\partial\overline\partial \delta\Phi$, we get
\begin{align}
\omega_{\Phi+\epsilon \delta\Phi}^d = \left(1 + \epsilon \Delta_{\overline{\partial}} \delta\Phi\right)\omega_\Phi^d + O(\epsilon^2) .
\end{align}
Expanding the logarithm with a Taylor expansion,
\begin{align}
\log\frac{\omega_{\Phi+\epsilon\,\delta\Phi}^d}{\omega^d} = \log\frac{\omega_\Phi^d}{\omega^d} + \epsilon \Delta_{\overline{\partial}}\delta\Phi + O(\epsilon^2) .
\end{align}
In particular, we have noted
\begin{align}
\frac{\omega_{\Phi+\epsilon\,\delta\Phi}^d}{\omega^d} = \frac{(1 + \epsilon \Delta_{\overline{\partial}} \delta\Phi)\omega_\Phi^d + O(\epsilon^2)}{\omega^d}  ,
\end{align}
and using log properties and Taylor expansions
\begin{align}
\log\frac{\omega_{\Phi+\epsilon\,\delta\Phi}^d}{\omega^d} = \log\left[ \frac{\omega_\Phi^d}{\omega^d} (1 + \epsilon \Delta_{\overline{\partial}} \delta\Phi) \right] & = \log\frac{\omega_\Phi^d}{\omega^d} + \log(1 + \epsilon \Delta_{\overline{\partial}} \delta\Phi) 
\\
& = \log\frac{\omega_\Phi^d}{\omega^d} + \epsilon \Delta_{\overline{\partial}}\delta\Phi + O(\epsilon^2) .
\end{align}
Now, let $u_\epsilon = \log\frac{\omega_{\Phi+\epsilon\,\delta\Phi}^d}{\omega^d} - f$, $V_\epsilon = \frac{\omega_{\Phi+\epsilon\,\delta\Phi}^d}{d!}$. At $\epsilon = 0$, we define
\begin{align}
\delta u = \Delta_{\overline{\partial}}\delta\Phi, \delta V = \frac{(\Delta_{\overline{\partial}}\delta\Phi)\omega_\Phi^d}{d!} .
\end{align}
The total variation obeys
\begin{align}
\delta A = \int_M (\delta u) V_0 + \int_M u_0 \delta V .
\end{align}
Substituting in the integral,
\begin{align}
\delta A = \int_M \left( \Delta_{\overline{\partial}}\delta\Phi \right) \frac{\omega_\Phi^d}{d!} + \int_M \left(\log\frac{\omega_\Phi^d}{\omega^d} - f\right) \left(\Delta_{\overline{\partial}}\delta\Phi\right) \frac{\omega_\Phi^d}{d!} .
\end{align}
Grouping,
\begin{align}
\delta A = \int_M \left( 1 + \log\frac{\omega_\Phi^d}{\omega^d} - f \right) \left( \Delta_{\overline{\partial}}\delta\Phi \right) \frac{\omega_\Phi^d}{d!} .
\end{align}
The first term vanishes since $\int \Delta_{\overline{\partial}} \delta {\Phi} = 0$ on a closed manifold due to the divergence theorem, and the second term with Green's second identity
\begin{align}
\delta A = \int_M \delta {\Phi} \cdot \Delta_{\overline{\partial}} \left( \log \frac{\omega_\Phi^d}{\omega^d} - f \right) \frac{\omega_\Phi^d}{d!}  = \int_M \left( \log\frac{\omega_\Phi^d}{\omega^d} - f \right) (\Delta_{\overline{\partial}}\delta\Phi) \frac{\omega_\Phi^d}{d!} = - \left\langle \log\frac{\omega_\Phi^d}{\omega^d} - f, \delta\Phi \right\rangle_\Phi .
\end{align}
Therefore, using the definition of the Dirichlet metric and the Riesz formulation we saw earlier, we establish, also using the definition of a gradient flow
\begin{align}
\dot{\Phi} = -\doublenabla \mathcal{F}_K = \log \frac{\omega_\Phi^d}{\omega^d} - f  .
\end{align}
Here, we have used the Dirichlet metric
\begin{align}
\langle \xi, \eta \rangle_{\Phi} = -\int_M \xi \cdot \Delta_{\overline{\partial}} \eta \frac{ \omega_{\Phi}^d }{d!} ,
\end{align}
which has a negative in its definition after integration by parts, therefore our result implicitly uses a double negative. The forms have the relationship
\begin{align}
\label{eqn:perelman_form_relationship}
\frac{\omega_\Phi^d}{\omega^d} = \frac{\det(\partial_i \partial_{\overline{j}} \Phi)}{\det(h_{0,i\overline{j}})} = \det \left(h^{0, i\overline{k}} \partial_{j} \partial_{\overline{k}} \Phi \right) =  \det (\partial_{i} \partial_{\overline{j}} \Phi ) ,
\end{align}
where the last equality follows when $h_0$ corresponds to a unit complex Gaussian. Taking the log,
\begin{align}
\log \left( \frac{\omega_\Phi^d}{\omega^d} \right) = \log \det \partial_{i} \partial_{\overline{j}} \Phi = \log \det h,
\end{align}
which is exactly used in the normalizing flow. Therefore,
\begin{align}
\dot{\Phi} = -\doublenabla \mathcal{F}_K =  \log \det h - f    .
\end{align}

\vspace{2mm}

\noindent Note that similar results that discuss
\begin{align}
\log \det h - f, \ \ \ \ \ \log \frac{\omega_\Phi^d}{\omega^d} - f 
\end{align}
in relation to Kähler potentials can be found in \cite{Cao1985} \cite{phong2007multiplieridealsheaveskahlerricci} \cite{Huang_2020} \cite{chen2009stabilitykahlerricciflow}, and various other references. The results in \cite{Cao1985} \cite{phong2007multiplieridealsheaveskahlerricci} \cite{Huang_2020} are mostly for an alternative smooth correction potential (not the Kähler potential). It can be noted that \cite{Cao1985} \cite{Huang_2020} do not mention our choice of metric. \cite{chen2009stabilitykahlerricciflow} does indeed present the result for the Kähler potential but also for the normalized Kähler-Ricci flow, and nor does it mention the Dirichlet metric. Most similar references are for the normalized Kähler-Ricci flow. 

\vspace{2mm}

\noindent Dirichlet metrics in the contexts of Kähler manifolds is primarily discussed in \cite{Calamai_2015}. \cite{Calamai_2015} also discusses Mabuchi and Calabi metrics, and defines the Dirichlet metric as
\begin{align}
\int_M \langle d\psi, d\chi\rangle_{\omega_t} \frac{\omega_{\Phi}^d}{d!} ,
\end{align}
which is equivalent to ours via integration by parts with vanishing boundary. \cite{Calamai_2015} does indeed discuss gradient flows under the Dirichlet metric and refers to \cite{Chen_2013}, which mostly discusses Calabi flows under functionals different from ours. This section concludes Theorem 2.

\noindent $ \square $

\section{Surgery of normalizing flows}

In this section, we describe the application of (modified, pseudo-) surgery as in Ricci flow \cite{perelman2003ricciflowsurgerythreemanifolds} in the context of complex normalizing flows. This section may have applications because typically singularity is a failure mode: the model is no longer bijective, and the model may crash. We consider a scenario when the normalizing flow diffeomorphism property collapses when the metric collapses, i.e. the Jacobian vanishes. 

\vspace{2mm}

\noindent We are interested in the case of singularity, i.e. collapse of the information metric
\begin{align}
\det(I_{i\overline{j}}(z, t)) \to 0 \quad \text{as} \quad t \rightarrow t^* .
\end{align}
In our contexts, this means that
\begin{align}
\det(h_{k^*})  \rightarrow 0 .
\end{align}
We also have that surgery can be triggered under the global condition, with local coordinates,
\begin{align}
\lim_{t \to t^*} \int_{M} |\text{Ric}(h)|  \frac{\omega_t^d}{d!} = +\infty  ,
\end{align} 
although sometimes this condition is local, which is a neck singularity where the normalizing flow ceases to be a biholomorphism. Now, at time $t^*$, we replace the degenerate metric with a new metric $\widetilde{h}_{i \overline{j}}$ defined via a Kähler potential perturbation $\varphi$
\begin{align}
\widetilde{h}_{i\overline{j}} = h_{i\overline{j}} + \partial_i \partial_{\overline{j}} \varphi .
\end{align}
$\varphi$ is smooth and $\det (\widetilde{h}) > 0 $. The surgery corresponds to a density transformation
\begin{align}
\log q_{k^*,\theta}(z) - \log q_{0,\theta}(z) \quad \mapsto  \quad \log q_{k^*,\theta}(z) - \log q_{0,\theta}(z) + \varphi(z) .
\end{align}
Based on our theory as in Section \ref{sec:main_contribution}, we desire across Ricci curvature accumulation
\begin{align}
\partial_i \partial_{\overline{j}} \varphi > - \Bigg( h_{0,i\overline{j}} + \mathbb{E}_{\theta} [\text{Ric}_{i\overline{j}}(h_{k^*}) - \text{Ric}_{i\overline{j}}(h_{0})] \Bigg),
\end{align}
which yields a complex Monge-Ampère inequality for the validity of the surgery
\begin{align}
\det\left(h_{k^*,i \overline{j}} + \partial_i \partial_{\overline{j}} \varphi\right) > 0 ,
\end{align}
and the flow restarts as
\begin{align}
 h_{\text{new},i\overline{j}} = h_{0,i\overline{j}} + \mathbb{E}_{\theta} [\text{Ric}_{i\overline{j}}(h_{k^*}) - \text{Ric}_{i\overline{j}}(h_{0})] + \partial_i \partial_{\overline{j}} \varphi  .
\end{align}
Let us consider the exponential of the log, which prevents the framework from being multiplicative. If it is multiplicative, the vanishing of the Jacobian map will still force the entire density to vanish, but if it additive, the potential now contributes nontrivially. In practice, we can compute easily
\begin{align}
\mathcal{L}_{k^*}(z) = \log q_{0,\theta}(z) - \sum_{k=1}^{k^*} \log |\det \mathcal{J}_k| ,
\end{align}
but this will negatively diverge at singularity, i.e.
\begin{align}
\lim_{|J_k| \rightarrow 0^+} \mathcal{L}_{k^*}(z)  = - \infty.    
\end{align}
We clip this, and at singularity, we consider the log-sum exponential
\begin{align}
\widetilde{\mathcal{L}}_{k^*}(z) = \exp \left\{ \text{LSE}(\mathcal{L}_{k^*}(z), \varphi(z)) \right\} =  \exp \mathcal{L}_{k^*}(z) + \exp \varphi(z)   ,
\end{align}
yielding the density
\begin{align}
\widetilde{q}_{k^*,\theta}(z) & = \frac{  \exp \mathcal{L}_{k^*}(z) + \exp \varphi(z)  }{  \int_{M} \left( e^{\mathcal{L}_{k^*}(z)} + e^{\varphi(z)} \right)  \det h \left( \frac{\sqrt{-1}}{2} \right)^d dz^1 \wedge d\overline{z}^1 \wedge \dots \wedge dz^d \wedge d\overline{z}^d }  .
\end{align}
Therefore, the inclusion of $\varphi$ counteracts the Jacobian collapse. Thus, at singularity, we recover up to a clipping
\begin{align}
\lim_{|J_k| \rightarrow 0^+} & \frac{  \exp \mathcal{L}_{k^*}(z) + \exp \varphi(z)  }{ \int_{M}  \left(  e^{\mathcal{L}_{k^*}(z)} + e^{\varphi(z)} \right) \det h \left( \frac{\sqrt{-1}}{2} \right)^d \bigwedge_i dz^i \wedge d\overline{z}^i }    = \frac{   \exp \varphi(z)  }{ \int_{M}  \exp \varphi(z)   \det h\left( \frac{\sqrt{-1}}{2} \right)^d \bigwedge_i dz^i \wedge d\overline{z}^i } .
\end{align}
We have noted $\det h$ on a set of nonzero measure does not imply the integral is zero over the entire domain, thus the integrals do not vanish even though $\det h $ can vanish locally. We remark the above application of the limit is also slightly informal, omitting limit properties and exchange for simplicity.

\section{Additional figures}

\begin{figure}[H]
  \vspace{0mm}
  \centering
  \includegraphics[scale=0.65]{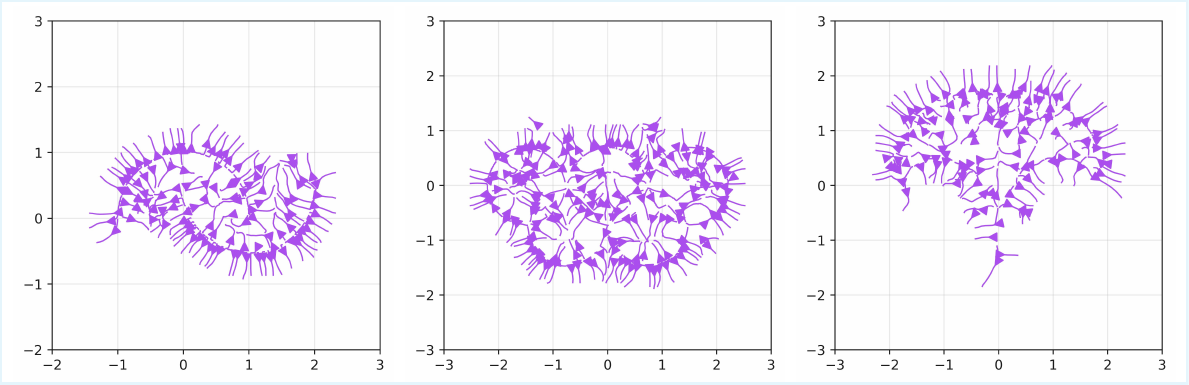}
  \caption{We plot $(\partial_x \Phi, \partial_y \Phi)$ at $t=T$ on $\mathbb{R}^2$, where $x = \text{Re}(z), y = \text{Im}(z)$, and $\Phi_t = -\log q_t$ (the Boltzmann constant $Z$ is omitted) using the baseline complex normalizing flow. Streamlines show the drift of the Langevin dynamics.}
  \label{fig:phi_gradients}
\end{figure}

\begin{figure}[H]
  \centering
  \vspace{10mm}
  \makebox[\textwidth][c]{\includegraphics[width=0.9\textwidth, page=1]{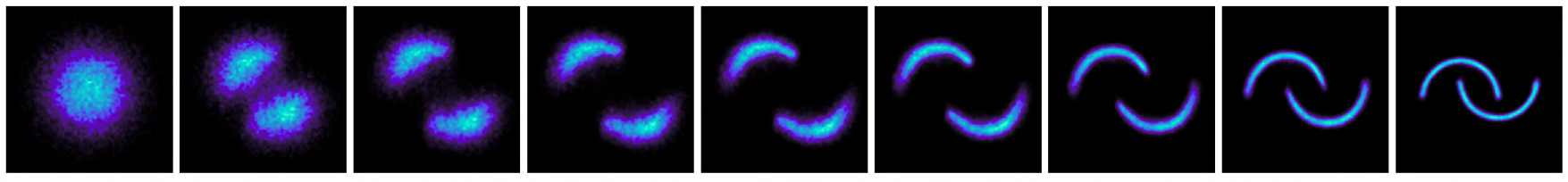}}\\[0.1em]
  \makebox[\textwidth][c]{\includegraphics[width=0.9\textwidth, page=2]{complexnf_layers_continuous.pdf}}\\[0.1em]
  \makebox[\textwidth][c]{\includegraphics[width=0.9\textwidth, page=3]{complexnf_layers_continuous.pdf}}
  \vspace{-2mm}
  \caption{We plot the output $(\text{Re}(z_k),\text{Im}(z_k))$ per $\Psi_{k,\theta}$ on our three datasets using a complex continuous normalizing flow.}
\end{figure}

\vspace{10mm}

\begin{figure}[H]
  \makebox[\textwidth][c]{\includegraphics[width=0.9\textwidth, page=1]{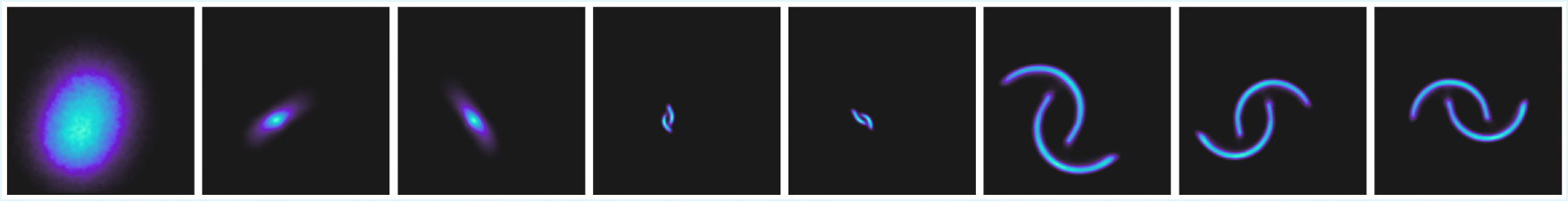}}\\[0.1em]
  \makebox[\textwidth][c]{\includegraphics[width=0.9\textwidth, page=2]{complexnf_layers_baseline.pdf}}\\[0.1em]
  \makebox[\textwidth][c]{\includegraphics[width=0.9\textwidth, page=3]{complexnf_layers_baseline.pdf}}
  \vspace{-2mm}
  \caption{We plot the output $(\text{Re}(z_k),\text{Im}(z_k))$ per $\Psi_{k,\theta}$ on our three datasets using a baseline complex normalizing flow. Empirically, we found the discrete flow worked better on the fractal tree (for continuous, we experimented with FFJORD \cite{grathwohl2018ffjordfreeformcontinuousdynamics} and custom; it is generally a known phenomenon that continuous normalizing flows are not particularly strong for the fractal tree, which is observable in \cite{zhang2021diffusionnormalizingflow}, although we remark \cite{zhang2021diffusionnormalizingflow} works well for fractal tree, which is also a diffusion model).}
  \label{fig:complexnf_layers_comparison}
\end{figure}

\end{document}